%% file: jel.tex
\overfullrule=0pt
\let\printlabel=N
\let\printdate=N
\tolerance=500
\input xpehfnt
\input xpehcmd
\let\afterstate=Y
\def\thestate{\the\sectionnum.\the\statenum}
\absskip=3mm
\def\BLF{\mathop{\hbox{\rm GLF}}}
\def\BBLF{\mathop{\hbox{\rm BGLF}}}
\def\GLF{GL-function}
\def\GLFs{GL-functions}
\def\BGLF{BGL-function}
\def\BGLFs{BGL-functions}
\def\UGLF{UGL-function}
\def\UGLFs{UGL-functions}
\def\PGLF{PGL-function}
\def\PGLFs{PGL-functions}
\def\GLS{GL-sequence}
\def\GLSs{GL-sequences}
\def\BGLS{BGL-sequence}

\def\GLFa{GL-family}
\def\GLFas{GL-families}
\def\Eig{\mathop{\hbox{\rm Eig}}}
\def\Dlim{\mathop{\hbox{\rm D-lim}}}
\def\Clim{\mathop{\hbox{\rm C-lim}}}
\def\Dlimsup{\mathop{\hbox{\rm D-limsup}}}
\def\Climsup{\mathop{\hbox{\rm C-limsup}}}
\def\D{\mathop{\hbox{\rm d}}}
\def\sD{\mathop{\hbox{\fnt{rm7}d}}}
\def\Du{\mathop{\hbox{\rm d}^{*}}}

\def\Dl{\mathop{\hbox{\rm d}_{*}}}

\def\range{\mathop{\hbox{\rm range}}}
\def\gcd{\mathop{\hbox{\rm g.c.d.}}}
\def\tf{\tilde{f}}
\def\af{\bar{f}}
\def\atf{\ovr{\tf}}
\def\ag{\bar{g}}
\def\hf{\hat{f}}
\def\tF{\widetilde{F}}
\def\hF{\widehat{F}}
\def\hphi{\hat{\phi}}
\def\tP{\widehat{P}}
\def\tv{\hat{v}}
\def\tw{\hat{w}}
\def\B{{\cal B}}
\def\T{{\cal T}}
\def\S{{\cal S}}
\def\gT{\mathord{\hbox{\fnt{gt}T}}}
\def\gTR{\gT_{\R}}
\def\TT{{\mb T}}
\def\H{{\cal H}}
\let\Lam=\Lambda
\def\Po{\mathord{\mathop{P}\limits^{\vbox{\hbox{\fnt{rm7}o}\kern-1pt}}}}
\def\Hc{\H_{\hbox{\fnt{rm7}c}}}
\def\Hwm{\H_{\hbox{\fnt{rm7}wm}}}
\let\cds=\cdots
\def\app{\mathop{\approx}\limits}
\def\blan{\bigl\lan}
\def\bran{\bigr\ran}
\def\Blan{\Bigl\lan}
\def\Bran{\Bigr\ran}
\def\P{{\cal P}}

\def\cW{{\cal W}}
\def\cL{{\cal L}}
\def\cM{{\cal M}}
\def\cN{{\cal N}}
\def\cG{{\cal G}}
\let\ophi=\oldphi
\def\tT{\widetilde{\T}}

\def\ff{\tilde{f}}
\def\ae{a.e.\nasp\ }
\start
\title{Joint ergodicity along generalized linear functions} 
\author{V. Bergelson, A. Leibman, and Y. Son}
\support{Bergelson and Leibman were supported by NSF grant DMS-1162073.}
\titledate
\address{Vitaly Bergelson}{Department of Mathematics, The Ohio State University, 
Columbus, OH 43210, USA}{bergelson.1@osu.edu}
\address{Alexander Leibman}{Department of Mathematics, The Ohio State University, 
Columbus, OH 43210, USA}{leibman.1@osu.edu}
\address{Younghwan Son}{Faculty of Mathematics and Computer Science,
The Weizmann Institute of Science, 234 Herzl Street, Rehovot 7610001 Israel}%
{younghwan.son@weizmann.ac.il}
\abstract{A criterion of joint ergodicity 
of several sequences of transformations of a probability measure space $X$ 
of the form $T_{i}^{\phi_{i}(n)}$
is given for the case where $T_{i}$ are commuting measure preserving transformations of $X$
and $\phi_{i}$ are integer valued generalized linear functions,
that is, the functions formed from conventional linear functions
by an iterated use of addition, multiplication by constants, 
and the greatest integer function.
We also establish a similar criterion for joint ergodicity
of families of transformations depending of a continuous parameter,
as well as a condition of joint ergodicity of sequences $T_{i}^{\phi_{i}(n)}$ along primes.}
\section{S-Intro}{Introduction}
Let $(X,\B,\mu)$ be a probability measure space.
A measure preserving transformation $T\col X\ra X$ is said to be {\it weakly mixing\/}
if the transformation $T\times T$,
acting on the Cartesian square $X\times X$,
is ergodic.
The notion of weak mixing was introduced in \brfr{vN-K}
(for measure preserving flows)
and has numerous equivalent forms
(see, for example, \brfr{B-Ros} and \brfr{B-Gor}.)
The following result involving weak mixing
plays a critical role in Furstenberg's proof (\brfr{F-Sz}) of ergodic Szemer\'{e}di theorem
and forms a natural starting point for numerous further developments
(see \brfr{B-pet}, \brfr{psz}, \brfr{BM-PolySz}, \brfr{BHa}):
\theorem{P-wmA}{}{}
If $T$ is an invertible weakly mixing measure preserving transformation of $X$,
then for any $k\in\N$ and any $A_{0},A_{1},\ld,A_{k}\in\B$ one has
$$
\lim_{N\ra\infty}\frac{1}{N}\sum_{n=1}^{N} 
\mu\bigl(A_{0}\cap T^{-n}A_{1}\cap\cds\cap T^{-kn}A_{k}\bigr)=\prod_{i=0}^{k}\mu(A_{i}).
$$ 
\endtheorem

It is not hard to show that \rfr{P-wmA} has the following functional form.
(In accordance with the well established tradition we write $Tf$ for the function $f(Tx)$.)

\theorem{P-wmF}{}{}
If $T$ is an invertible weakly mixing measure preserving transformation of $X$,
then for any $k\in\N$, any distinct nonzero integers $a_{1},\ld,a_{k}$,
and any $f_{1},\ld,f_{k}\in L^{\infty}(X)$ one has
$$
\lim_{N\ra\infty}\frac{1}{N}\sum_{n=1}^{N}
T^{a_{1}n}f_{1}\cds T^{a_{k}n}f_{k}=\prod_{i=1}^{k}\int_{X}f_{i}\,d\mu
$$ 
in $L^{2}$ norm.
\endtheorem
\npar
In other words, given a weakly mixing transformation $T$ of $X$
and distinct nonzero integers $a_{1},\ld,a_{k}$,
the transformations $T^{a_{1}},\ld,T^{a_{k}}$
(or, rather, the sequences $T^{a_{1}n},\ld,T^{a_{k}n}$, $n\in\N$)
possess a strong independence property.
This naturally leads to the following definition:
\definition{} (Cf.\nasp\ \brfr{BB1}.)
Measure preserving transformations $T_{1},\ld,T_{k}$ of a probability measure space $X$
are said to be {\it jointly ergodic\/}
if for any $f_{1},\ld,f_{k}\in L^{\infty}(X)$ one has
$$
\lim_{N\ra\infty}\frac{1}{N}\sum_{n=1}^{N}
T_{1}^{n}f_{1}\cds T_{k}^{n}f_{k}=\prod_{i=1}^{k}\int_{X}f_{i}\,d\mu
$$ 
in $L^{2}$ norm.
\enddefinition

The following theorem, proved in \brfr{BB1},
provides a criterion of joint ergodicity of commuting measure preserving transformations:
\theorem{P-BB1}{}{}
Let $T_{1},\ld,T_{k}$ be commuting invertible measure preserving transformations of $X$.
Then $T_{1},\ld,T_{k}$ are jointly ergodic 
iff the transformation $T_{1}\times\cds\times T_{k}$ of $X^{k}$ is ergodic
and the transformations $T_{i}^{-1}T_{j}$ of $X$ are ergodic for all $i\neq j$.
\endtheorem

Further developments
(most of which were motivated by connections with combinatorics and number theory)
have revealed that the phenomenon of joint ergodicity is a rather general one.
For example, as it was shown in \brfr{B-pet},
if $T$ is an invertible weakly mixing measure preserving transformation
and $p_{1},\ld,p_{k}$ are nonconstant polynomials $\Z\ra\Z$ 
with $p_{i}-p_{j}\neq\const$ for any $i\neq j$,
then for any $f_{1},\ld,f_{k}\in L^{\infty}(X)$ 
one has
$$
\lim_{N\ra\infty}\frac{1}{N}\sum_{n=1}^{N}
T^{p_{1}(n)}f_{1}\cds T^{p_{k}(n)}f_{k}=\prod_{i=1}^{k}\int_{X}f_{i}\,d\mu
$$ 
in $L^{2}$ norm.
(See also \brfr{FraKra}, \brfr{BHa}, and \brfr{Fra} for more results of this flavor.)
So, it makes sense to consider ergodicity and joint ergodicity 
of sequences of measure preserving transformations of general form:
\definition{}
Let $\T(n)$, $n\in\N$, be a sequence of measure preserving transformations of $X$;
we say that $\T$ is {\it ergodic\/}
if for any $f\in L^{2}(X)$, 
$$
\lim_{N\ra\infty}\frac{1}{N}\sum_{n=1}^{N}\T(n)f\dsc=\int_{X}f\,d\mu.
$$
Given several sequences $\T_{1}(n),\ld,\T_{k}(n)$, $n\in\N$, 
of measure preserving transformations of $X$,
we say that $\T_{1},\ld,\T_{k}$ are {\it jointly ergodic\/} if
$$
\lim_{N\ra\infty}\frac{1}{N}\sum_{n=1}^{N} 
\T_{1}(n)f_{1}\cds\T_{k}(n)f_{k}=\prod_{i=1}^{k}\int_{X}f_{i}\,d\mu
$$
in $L^{2}$~norm for any $f_{1},\ld,f_{k}\in L^{\infty}(X)$.
\enddefinition

\ignore
Throughout the paper $(X,\B,\mu)$ is a probability measure space.
 
Measure preserving transformations $T_{1},\ld,T_{k}$ of $X$
are said to be {\it jointly ergodic\/} if
$$
\lim_{N\ra\infty}\frac{1}{N}\sum_{n=1}^{N} 
T_{1}^{n}f_{1}\cds T_{k}^{n}f_{k}=\prod_{i=1}^{k}\int_{X}f_{i}\,d\mu
$$ 
in the $L^2$~norm for any $f_{1},\ld,f_{k}\in L^{\infty}(X)$.
The following theorem, established in \brfr{BB1},
provides a criterion of joint ergodicity 
of several commuting measure preserving transformations:
\theorem{P-BB1}{}{}
Let $T_{1},\ld,T_{k}$ be commuting invertible measure preserving transformations of $X$.
Then $T_{1},\ld,T_{k}$ are jointly ergodic 
iff the transformation $T_{1}\times\cds\times T_{k}$ of $X^{k}$ is ergodic
and the transformations $T_{i}^{-1}T_{j}$ of $X$ are ergodic for all $i\neq j$.
\endtheorem
Let now $\T(n)$, $n\in\N$, be a sequence of measure preserving transformations of $X$;
we say that $\T$ is {\it ergodic\/}
if for any $f\in L^{2}(X)$, 
$$
\lim_{N\ra\infty}\frac{1}{N}\sum_{n=1}^{N}\T(n)f\dsc=\int_{X}f\,d\mu.
$$
Given several sequences $\T_{1}(n),\ld,\T_{k}(n)$, $n\in\N$, 
of measure preserving transformations of $X$,
we say that $\T_{1},\ld,\T_{k}$ are {\it jointly ergodic\/} if
$$
\lim_{N\ra\infty}\frac{1}{N}\sum_{n=1}^{N} 
\T_{1}(n)f_{1}\cds\T_{k}(n)f_{k}=\prod_{i=1}^{k}\int_{X}f_{i}\,d\mu
$$
in $L^{2}$~norm for any $f_{1},\ld,f_{k}\in L^{\infty}(X)$.
\endignore

Results obtained in \brfr{B-pet}, \brfr{BHa}, and \brfr{Fra}
lead to a natural question of what are the necessary and sufficient conditions
for joint ergodicity of sequences of transformations
of the form $T_{1}^{\phi_{1}(n)},\ld,T_{k}^{\phi_{k}(n)}$,
where $T_{i}$ are measure preserving transformations of $X$
and $\phi_{i}(n)$ are ``sufficiently regular'' sequences of integers diverging to infinity.
In the case where $T_{1}=\ld=T_{k}=T$ where $T$ is a weakly mixing transformation, 
this question has a quite satisfactory answer 
not only when $\phi_{i}$ are integer-valued polynomials,
but also, more generally, are functions of the form $[\psi_{i}]$,
where $[\cd]$ denotes the integer part
and $\psi_{i}$ are either the so-called ``tempered functions'',
or functions of polynomial growth belonging to a Hardy field 
(see \brfr{BHa} and \brfr{Fra}).

Much less is known about joint ergodicity of $T_{1}^{\phi_{1}(n)},\ld,T_{k}^{\phi_{k}(n)}$
when $T_{i}$ are distinct, not necessarily weakly mixing transformations.
It is our goal in this paper to extend \rfr{P-BB1}
to the case $\phi_{i}$ are integer-valued {\it generalized linear functions}.
A generalized, or {\it bracket\/} linear function (of real or integer argument)
is a function constructible from conventional linear functions
with the help of the operations of addition, multiplication by constants,
and taking the integer part, $[\cd]$ (or, equivalently, the fractional part, $\{\cd\}$).
(For example, $\phi(n)=[\alf_{1}n+\alf_{2}]$, 
$\phi(n)=\alf_{1}[\alf_{2}n+\alf_{3}]+\alf_{4}$,
and, say, $\phi(n)=\alf_{1}\bigl[\alf_{2}\bigl[\alf_{3}[\alf_{4}n+\alf_{5}]+\alf_{6}\bigr]
+\alf_{7}[\alf_{8}n+\alf_{9}]\bigr]+\alf_{10}n+\alf_{11}$, 
where $\alf_{i}\in\R$,
are generalized linear functions.)
In complete analogy with \rfr{P-BB1}, we prove:
\theorem{P-reg}{}{}
Let $T_{1},\ld,T_{k}$ be commuting invertible measure preserving transformations of $X$
and let $\phi_{1},\ld,\phi_{k}$ be generalized linear functions $\Z\ra\Z$.
The sequences $T_{1}^{\phi_{1}(n)},\ld,T_{k}^{\phi_{k}(n)}$ are jointly ergodic 
iff the sequence $T_{1}^{\phi_{1}(n)}\times\cds\times T_{k}^{\phi_{k}(n)}$ 
of transformations of $X^{k}$ is ergodic
and the sequences $T_{i}^{-\phi_{i}(n)}T_{j}^{\phi_{j}(n)}$ of transformations of $X$ 
are ergodic for all $i\neq j$.
\endtheorem
Here are two special cases of \rfr{P-reg}:
\corollary{P-single}{}{}
Let $T$ be a weakly mixing invertible measure preserving transformation of $X$
and let $\phi_{1},\ld,\phi_{k}$ be unbounded generalized linear functions $\Z\ra\Z$.
such that $\phi_{j}-\phi_{i}$ are unbounded for all $i\neq j$. 
Then for any $f_{1},\ld,f_{k}\in L^{\infty}(X)$,
$$
\lim_{N\ras\infty}\frac{1}{N}\sum_{n=1}^{N}T^{\phi_{1}(n)}f_{1}\cds T^{\phi_{k}(n)}f_{k}
=\prod_{i=1}^{k}\int_{X}f_{i}\,d\mu.
$$
In particular, for any distinct $\alf_{1},\ld,\alf_{k}\in\R\sm\{0\}$,
$$
\lim_{N\ras\infty}\frac{1}{N}\sum_{n=1}^{N}T^{[\alf_{1}n]}f_{1}\cds T^{[\alf_{k}n]}f_{k}
=\prod_{i=1}^{k}\int f_{i}\,d\mu.
$$
\endcorollary
For a measure preserving transformation $T$ of $X$,
let $\Eig T$ be the set of eigenvalues of $T$,
$$
\Eig T=\bigl\{\lam\in\C^{*}:\hbox{$Tf = \lam f$ for some $f\in L^{2}(X)$}\bigr\}.
$$
For several measure preserving transformations $T_{1},\ld,T_{k}$ of $X$
we put $\Eig(T_{1},\ld,T_{k})=\prod_{i=1}^{k}\Eig T_{i}$.

\corollary{P-spec}{}{}
Let $T_{1},\ld,T_{k}$ be commuting invertible 
jointly ergodic measure preserving transformations of $X$
and let $\phi$ be an unbounded generalized linear function $\Z\ra\Z$.
Then
$$
\lim_{N\ras\infty}\frac{1}{N}\sum_{n=1}^{N}T_{1}^{\phi(n)}f_{1}\cds T_{k}^{\phi(n)}f_{k}
=\prod_{i=1}^{k}\int f_{i}\,d\mu\
\hbox{for any $f_{1},\ld,f_{k}\in L^{\infty}(X)$}
$$
iff $\lim_{N\ras\infty}\frac{1}{N}\sum_{n=1}^{N}\lam^{\phi(n)}=0$
for every $\lam\in\Eig(T_{1},\ld,T_{k})\sm\{1\}$.
In particular, for any irrational $\alf\in\R$,
$$
\lim_{N\ras\infty}\frac{1}{N}\sum_{n=1}^{N}T_{1}^{[\alf n]}f_{1}\cds T_{k}^{[\alf n]}f_{k}
=\prod_{i=1}^{k}\int f_{i}\,d\mu\
\hbox{for any $f_{1},\ld,f_{k}\in L^{\infty}(X)$}
$$ 
iff $e^{2\pi i\alf^{-1}\Q}\cap\Eig(T_{1},\ld,T_{k})=\{1\}$.
\endcorollary
In fact, we obtain a result more general than \rfr{P-reg}.
Let $G$ be a commutative group of measure preserving transformations of $X$.
We say that a sequence $\T$ of transformations of $X$
is a {\it generalized linear sequence\/} in $G$
if it has the form $\T(n)=T_{1}^{\phi_{1}(n)}\cds T_{r}^{\phi_{r}(n)}$, $n\in\Z$,
for some $T_{1},\ld,T_{r}\in G$ 
and generalized linear functions $\phi_{1},\ld,\phi_{k}\col\Z\ra\Z$.
(The sequences $T_{i}^{-\phi_{i}(n)}T_{j}^{\phi_{j}(n)}$ appearing in \rfr{P-reg} 
are of this sort.)
Also, we change the definitions of ergodicity and of joint ergodicity above,
replacing the averages $\frac{1}{N}\sum_{n=1}^{N}$
with the more general averages $\frac{1}{|\Phi_{N}|}\sum_{n\in\Phi_{N}}$,
where $(\Phi_{N})$ is an arbitrary F{\o}lner sequence in $\Z$.
(See \rfr{D-jerg}.)
The {\it uniform\/} ergodicity and joint ergodicity,
which appear when the averages $\frac{1}{N}\sum_{n=1}^{N}$, with $N\ras\infty$,
are replaced by the averages $\frac{1}{M-N}\sum_{n=N+1}^{M}$, with $M-N\ras\infty$,
form a special case of it.)
In this setup, we prove the following:
\theorem{P-main}{}{}
Generalized linear sequences $\T_{1},\ld,\T_{k}$ 
in a commutative group of transformations of $X$ are jointly ergodic
iff the sequence $\T_{1}\times\cds\times\T_{k}$ of transformations of $X^{k}$ is ergodic
and the sequences $\T_{i}^{-1}\T_{j}$ of transformations of $X$ 
are ergodic for all $i\neq j$.
\endtheorem
In addition to \rfr{P-reg}, we also prove a version thereof along primes.
In particular, we obtain the following result:
\theorem{P-reg-primes}{}{}
Let $T_{1},\ld,T_{k}$ be commuting invertible measure preserving transformations of $X$
and let $\phi_{1},\ld,\phi_{k}$ be generalized linear functions $\Z\ra\Z$.
Assume that for any $W\in\N$ and $r\in R(W)$
the sequences $T_{i}^{\phi_{i}(Wn+r)}$, $i=1,\ld,k$, are jointly ergodic.
Then for any $f_{1},\ld,f_{k}\in L^{\infty}(X)$,
$$
\lim_{N\ras\infty}\frac{1}{\pi(N)}\sum_{p\in\P(N)}
T_{1}^{\phi_{1}(p)}f_{1}\cds\T_{k}^{\phi_{k}(p)}f_{k}=\prod_{i=1}^{k}\int_{X}f_{i}\,d\mu
\quad\hbox{(in $L^{2}$ norm)}.
$$
\endtheorem
The structure of the paper is as follows:
Sections \rfrn{S-GLF}-\rfrn{S-APFB} contain technical material 
related to properties of generalized linear functions.
In \rfr{S-GLST} we investigate ergodic properties 
of what we call ``a generalized linear sequence of measure preserving transformations''
-- a product of several sequences of the form $T^{\phi(n)}$,
where $\phi$ is an integer valued generalized linear function.
In \rfr{S-JE} we obtain our main result, \rfr{P-joint},
the criterion of joint ergodicty of several commuting generalized linear sequences.
In \rfr{S-Primes}, we extend \rfr{P-joint} to averaging along primes.
In \rfr{S-cont} we deal with families of transformations depending on a continuous parameter,
and obtain a version of \rfr{P-joint} for continuous flows.
By using a ``change of variable'' trick we also extend this result 
to more general families of transformations of the form $T^{\phi(\sig(t))}$, 
where $\phi$ is a generalized linear function
and $\sig$ is a monotone function of ``regular'' growth.
For example, we have the following version of \rfr{P-spec}:
\proposition{P-IRspec}{}{}
Let $T_{1}^{s},\ld,T_{k}^{s}$, $s\in\R$, 
be commuting jointly ergodic continuous flows of measure preserving transformations of $X$
and let $\phi$ be an unbounded generalized linear function;
then for any $\alf>0$ 
the families $T_{1}^{\phi(t^{\alf})},\ld,T_{k}^{\phi(t^{\alf})}$, $t\in[0,\infty)$, 
are jointly ergodic
(that is, $\lim_{b\ras\infty}\frac{1}{b}\int_{0}^{b}
T_{1}^{\phi(t^{c})}f_{1}\cds T_{k}^{\phi(t^{c})}f_{k}\,dt=\prod_{i=1}^{k}\int_{X}f_{i}\,d\mu$
for any $f_{1},\ld,f_{k}\in L^{\infty}(X)$)
iff $\lim_{b\ras\infty}\frac{1}{b}\int_{0}^{b}\lam^{\phi(t)}\,dt=0$ 
for every $\lam\in\Eig(T_{1}^{1},\ld,T_{k}^{1})\sm\{1\}$.
\endproposition
\npar
Finally, \rfr{S-noncom} contains a result 
pertaining to joint ergodicity of several non-commuting generalized linear sequences.
\section{S-GLF}{Generalized linear functions}

For $x\in\R$ we denote by $[x]$ the integer part of $x$
and by $\{x\}$ the fractional part $x-[x]$ of $x$.

The set $\BLF$ of {\it generalized linear functions\/}
is the minimal set of functions $\R\ra\R$ 
containing all linear functions $ax+b$
and closed under addition, multiplication by constants,
and the operation of taking the integer (equivalently, the fractional) part.
More exactly, we define $\BLF$ inductively in the following way.
We put $\BLF_{0}=\bigl\{\phi(x)=ax+b,\ a,b\in\R\bigr\}$.
After $\BLF_{k}$ has already been defined,
we define $\BLF_{k+1}$ to be the space of functions
spanned by $\BLF_{k}$ and the set $\bigl\{[\phi],\ \phi\in\BLF_{k}\bigr\}$.
(Equivalently, we can define $\BLF_{k+1}$ to be the space
spanned by $\BLF_{k}$ and the set $\bigl\{\{\phi\},\ \phi\in\BLF_{k}\bigr\}$.)
Finally, we put $\BLF=\bigcup_{k=0}^{\infty}\BLF_{k}$.
For $\phi\in\BLF$, we call the minimal $k$ for which $\phi\in\BLF_{k}$
{\it the weight\/} of $\phi$.

We will refer to functions from $\BLF$ as to {\it \GLFs}.
\example{}
$\phi(x)=a_{1}\bigl\{a_{2}\bigl[a_{3}\{a_{4}x+a_{5}\}+a_{6}\bigr]
+a_{7}[a_{8}x+a_{9}]\bigr\}+a_{10}x+a_{11}$, 
where $a_{1},\ld,a_{11}\in\R$, is a \GLF.
\endexample
Clearly, the set of \GLFs\ is closed under the composition:
if $\phi_{1},\phi_{2}\in\BLF$, then $\phi_{1}(\phi_{2}(x))\in\BLF$.

We define the set $\BBLF$ inductively in the following way:
$\BBLF_{1}=\bigl\{\phi(x)=\{ax+b\},\ a,b\in\R\bigr\}$;
if $\BBLF_{k}$ has already been defined,
$\BBLF_{k+1}$ is the space 
spanned by the set $\BBLF_{k}\cup\bigl\{\{\phi\},\ \phi\in\BLF_{k}\bigr\}$;
and finally, $\BBLF=\bigcup_{k=1}^{\infty}\BBLF_{k}$.

\lemma{P-def\BGLF}{}{}
$\BBLF$ is exactly the set of bounded \GLFs.
(Hence the abbreviation ``BGLF''.)
\endlemma
\proof{}
Clearly, all elements of $\BBLF$ are bounded \GLFs.
To prove the opposite inclusion
we use induction on the weight of \GLFs.
Let $\phi\in\BLF_{k}\sm\BLF_{k-1}$ be bounded.
If $k=0$, then $\phi$ must be a constant and thus belongs to $\BBLF$.
If $k\geq 1$, $\phi=\phi_{0}+\sum_{i=1}^{m}a_{i}\{\phi_{i}\}$,
where $\phi_{0},\phi_{1},\ld,\phi_{m}\in\BLF_{k-1}$.
Now, $\phi_{0}$ is bounded, thus by induction, $\phi_{0}\in\BBLF$,
and $\{\phi_{1}\},\ld,\{\phi_{m}\}\in\BBLF$ by definition,
so $\phi\in\BBLF$.%
\endproof
We will refer to elements of $\BBLF$ 
as to {\it bounded generalized linear functions}, or {\it \BGLFs}.
\lemma{P-LBP}{}{}
Any \GLF\ $\phi$ is uniquely representable in the form $\phi(x)=ax+\psi(x)$, 
where $a\in\R$ and $\psi$ is a \BGLF.
\endlemma
\proof{}
Every $\phi\in\BLF_{k}$ has the form $\phi=\phi_{0}+\sum_{i=1}^{m}a_{i}\{\phi_{i}\}$
with $\phi_{0},\phi_{1},\ld,\phi_{m}\in\BLF_{k-1}$.
We have $\sum_{i=1}^{m}a_{i}\{\phi_{i}\}\in\BBLF$,
and $\phi_{0}$ is representable in the form $\phi_{0}(x)=ax+\psi_{0}(x)$ with $\psi_{0}\in\BBLF$ 
by induction on $k$.

As for the uniqueness, if $a_{1}x+\psi_{1}(x)=a_{2}x+\psi_{2}(x)$
with $a_{1},a_{2}\in\R$ and $\psi_{1},\psi_{2}\in\BBLF$,
then the function $(a_{1}-a_{2})x$ is bounded, and so $a_{1}=a_{2}$.%
\endproof
\corollary{P-Intp}{}{}
Any \GLF\ $\phi$ is uniquely representable in the form $\phi(x)=[ax]+\xi(x)$
with $a\in\R$ and $\xi\in\BBLF$.
\endcorollary
For a function $\phi\col\R\ra\R$ and $\alf\in\R$
``the difference derivative'' $D_{\alf}\phi$ of $\phi$ with step $\alf$
is $D_{\alf}\phi(x)=\phi(x+\alf)-\phi(x)$, $x\in\R$.
\corollary{P-Dif}{}{}
For any \GLF\ $\phi$ and $\alf\in\R$, $D_{\alf}\phi$ is a \BGLF.
\endcorollary
We will refer to \BGLFs\ taking values in $\{0,1\}$ as to {\it \UGLFs}.
\lemma{P-ind}{}{}
Let $\phi$ be a \BGLF.
Then for any $a\in\R$
the indicator functions 
$1_{\{\phi<a\}}$, $1_{\{\phi\leq a\}}$, $1_{\{\phi>a\}}$, and $1_{\{\phi\geq a\}}$
of the sets
$\{x:\phi(x)<a\}$, $\{x:\phi(x)\leq a\}$, $\{x:\phi(x)>a\}$, and $\{x:\phi(x)\geq a\}$
are \UGLFs.
\endlemma
\proof{}
We start with the set $\{\phi\geq a\}$.
Let $c=\sup|\phi|+|a|+1$.
Then the function $\xi=(\phi-a)/c+1$ satisfies $0<\xi<2$,
and $\phi\geq a$ iff $\xi\geq 1$.
Thus, the \UGLF\ $[\xi]$ is just $1_{\{\phi\geq a\}}$.

Now, $1_{\{\phi\leq a\}}=1_{\{-\phi\geq-a\}}$,
$1_{\{\phi<a\}}=1-1_{\{\phi\leq a\}}$, and  $1_{\{\phi>a\}}=1-1_{\{\phi\leq a\}}$.%
\endproof
We will now show that the set of \UGLFs\ is closed under Boolean operations.
For two functions $\phi$ and $\psi$ taking values in $\{0,1\}$,
let $\phi\vee \psi=\max\{\phi,\psi\}=\phi+\psi-\phi\psi$, 
$\phi\wedge \psi=\min\{\phi,\psi\}=\phi\psi$, 
and $\neg \phi=1-\phi$.
\proposition{P-log}{}{}
If $\phi$, $\psi$ are \UGLFs, then $\phi\vee\psi$, $\phi\wedge\psi$, and $\neg\phi$ 
are also \UGLFs.
\endproposition
\proof{}
$\neg \phi=1-\phi$ is clearly a \UGLF,
$\phi\vee \psi$ is the indicator function of the set $\{\phi+\psi>0\}$
and thus is a \UGLF\ by \rfr{P-ind},
and $\phi\wedge \psi=\neg(\neg \phi\vee\neg \psi)$.%
\endproof
From \rfr{P-log} we get the following generalization of \rfr{P-ind}:
\proposition{P-indI}{}{}
Let $\phi_{1},\ld,\phi_{k}$ be \BGLFs\ and let $\phi=(\phi_{1},\ld,\phi_{k})$.
For any interval $I=I_{1}\times\cds\times I_{k}\sle\R^{k}$,
(where $I_{i}$ are intervals in $\R$, 
which may be bounded or unbounded, open, closed, half-open half-closed, or degenerate)
the indicator function $1_{A}$ of the set $A=\{x:\phi(x)\in I\}$ is a \UGLF.
\endproposition
We also have the following:
\proposition{P-char}{}{}
Let $\phi$ be an unbounded \GLF\ $\Z\ra\Z$.
Then the indicator function $1_{H}$ 
of the range $H=\phi(\Z)$ of $\phi$ is a \UGLF.
\endproposition
\npar(Notice that \GLFs\ $\Z\ra\R$ are restrictions of \GLFs\ $\R\ra\R$,
thus all the results above apply.)
\proof{}
By \rfr{P-Intp},
$\phi(n)=[an]+\psi(n)$ for some $a\in\R$, $\xi\in\BBLF$.
Since $\phi$ is integer-valued, $\psi$ is integer valued,
and thus the range $K=\psi(\Z)$ of $\psi$ is a finite set of integers.
Since $\phi$ is unbounded, $a\neq 0$;
let us assume that $a>0$.

If $n,k,j\in\Z$ are such that $n=[ak]+j$,
then $0\leq ak-n+j<1$, so
$$\txt
\frac{n-j}{a}\leq k<\frac{n-j+1}{a}
$$
and so, 
$$\txt
k\in\bigl\{\bigl[\frac{n-j}{a}\bigr]+i,\ i\in I\bigr\},
$$
where $I=\bigl\{0,1,\ld,\bigl[\frac{1}{a}\bigr]+1\bigr\}$.
Hence, if $n\in H$, that is, if $n=\phi(k)$ for some $k\in\Z$,
then 
$$\txt
k\in\bigl\{\bigl[\frac{n-j}{a}\bigr]+i,\ i\in I,\ j\in K\bigr\}.
$$

For each $i\in I$ and $j\in K$, define
$$\txt
\del_{i,j}(n)=n-\phi\bigl(\bigl[\frac{n-j}{a}\bigr]+i\bigr)
=n-a\bigl(\bigl[\frac{n-j}{a}\bigr]+i\bigr)
-\psi\bigl(\bigl[\frac{n-j}{a}\bigr]+i\bigr);
$$
then $\del_{i,j}\in\BBLF$ for all $i,j$,
and $n\in H$ iff $\del_{i,j}(n)=0$ for some $i,j$.
By \rfr{P-ind} and \rfr{P-log}
the indicator functions $1_{\{\del_{i,j}=0\}}$ are \UGLFs\ for all $i,j$,
and thus the function
$1_{H}=\bigvee_{\sdup{i\in I}{j\in K}}1_{\{\del_{i,j}=0\}}$
is also a \UGLF\ by \rfr{P-log}.%
\endproof
\section{S-Corp}{C-lims, D-lims, densities, and the van der Corput trick}

This is a technical section.
Starting from this moment
we fix an arbitrary F{\o}lner sequence $(\Phi_{N})_{N=1}^{\infty}$ in $\Z$
(that is, a sequence of finite subsets of $\Z$ with the property that for any $h\in\Z$,
$|(\Phi_{N}-h)\triangle\Phi_{N}|/|\Phi_{N}|\ras 0$ as $N\ras\infty$).

Under ``a sequence'' we will usually understand a function with domain $\Z$.
For a sequence $(u_{n})$ of real numbers, or of elements of a normed vector space,
we define 
$\Clim_{n}u_{n}=\lim_{N\ras\infty}\frac{1}{|\Phi_{N}|}\sum_{n\in\Phi_{N}}u_{n}$,
if this limit exists.
When $u_{n}$ are real numbers, we define
$\Climsup_{n}u_{n}=\limsup_{N\ras\infty}\frac{1}{|\Phi_{N}|}\sum_{n\in\Phi_{N}}u_{n}$.
When $u_{n}$ are elements of a normed vector space 
we also define
$\Climsup_{\|\cd\|,n}u_{n}
=\limsup_{N\ras\infty}\dsc\bigl\|\frac{1}{|\Phi_{N}|}\sum_{n\in\Phi_{N}}u_{n}\bigr\|$.

For a set $E\sle\Z$ we define {\it the density\/} of $E$ to be
$\D(E)=\lim_{N\ras\infty}|E\cap\Phi_{N}|/|\Phi_{N}|$,
if this limit exists.
We also define {\it the upper density} and {\it the lower density\/} of $E$
as $\Du(E)=\limsup_{N\ras\infty}|E\cap\Phi_{N}|/|\Phi_{N}|$
and $\Dl(E)=\liminf_{N\ras\infty}|E\cap\Phi_{N}|/|\Phi_{N}|$ respectively.

We will say that a sequence $(z_{n})$ in a probability measure space $(Z,\lam)$
is {\it uniformly distributed\/}
if $\Clim_{n}g(z_{n})=\int_{Z}g\,d\lam$ for any $g\in C(Z)$.

For a sequence $(u_{n})$ of vectors in a normed vector space
we write $\Dlim_{n}u_{n}=u$ if for any $\eps>0$,
$\D\bigl(\bigl\{n:\|u_{n}-u\|\geq\eps\bigr\}\bigr)=0$. 
Clearly, this is equivalent to $\Clim_{n}\|u_{n}-u\|=0$.
For a sequence $(u_{n})$ of real numbers we also define $\Dlimsup_{n}u_{n}$ 
as $\inf\bigl\{u\in\R:\D(\{n:u_{n}>u\})=0\bigr\}$.

We will be using the following version of the van der Corput trick:
\lemma{P-vdC}{}{}
Let $(u_{n})$ be a bounded sequence of elements of a Hilbert space.
Then for any finite subset $D$ of $\Z$,
\equ{
\Climsup\nolimits_{n,\|\cd\|}u_{n}
\leq\Bigl(\frac{1}{|D|^{2}}\sum_{h_{1},h_{2}\in D}
\Climsup_{n}\lan u_{n+h_{1}},u_{n+h_{2}}\ran\Bigr)^{1/2}.
}
Thus, if for some $\eps>0$ there exists an infinite set $B\sle\Z$ 
such that $\bigl|\Climsup_{n}\lan u_{n+h_{1}},u_{n+h_{2}}\ran\bigr|<\eps$ 
for all distinct $h_{1},h_{2}\in B$,
then $\Climsup_{n,\|\cd\|}u_{n}<\sqrt{\eps}$.
\endlemma
\proof{}
Let $D\sle\Z$, $|D|<\infty$.
For any $N\in\N$ we have
\equ{
\frac{1}{|\Phi_{N}|}\sum_{n\in\Phi_{N}}u_{n}
=\frac{1}{|D|}\sum_{h\in D}\frac{1}{|\Phi_{N}|}\sum_{n\in\Phi_{N}}u_{n}
=\Bigl(\frac{1}{|D|}\sum_{h\in D}
\frac{1}{|\Phi_{N}|}\sum_{n\in\Phi_{N}}u_{n+h}\Bigr)
-A_{N}+B_{N},
}
where $A_{N}=\frac{1}{|D|}\sum_{h\in D}
\frac{1}{|\Phi_{N}|}\sum_{\sdup{n\in\Phi_{N}}
{\hbox to 0pt{\hss$\scr n+h\not\in\Phi_{N}$\hss}}}u_{n+h}$
and $B_{N}=\frac{1}{|D|}\sum_{h\in D}
\frac{1}{|\Phi_{N}|}\sum_{\sdup{n\not\in\Phi_{N}}
{\hbox to 0pt{\hss$\scr n+h\in\Phi_{N}$\hss}}}u_{n+h}$.
Since $\{\Phi_{N}\}_{N=1}^{\infty}$ is a F{\o}lner sequence
and the sequence $(u_{n})$ is bounded,
$\|A_{N}\|,\|B_{N}\|\ras 0$ as $N\ras\infty$.
Thus, 
\equ{
\limsup_{N\ras\infty}\Bigl\|
\frac{1}{|\Phi_{N}|}\sum_{n\in\Phi_{N}}u_{n}\Bigr\|
=\limsup_{N\ras\infty}\Bigl\|\frac{1}{|D|}\sum_{h\in D}
\frac{1}{|\Phi_{N}|}\sum_{u\in\Phi_{N}}u_{n+h}\Bigr\|.
}
By Schwarz's inequality,
\lequ{
\Bigl\|\frac{1}{|D|}\sum_{h\in D}
\frac{1}{|\Phi_{N}|}\sum_{n\in\Phi_{N}}u_{n+h}\Bigr\|^{2}
=\frac{1}{|D|^{2}}\Bigl\|\frac{1}{|\Phi_{N}|}
\sum_{n\in\Phi_{N}}\sum_{h\in D}u_{n+h}\Bigr\|^{2}
\leq\frac{1}{|D|^{2}}\frac{1}{|\Phi_{N}|}
\sum_{n\in\Phi_{N}}\Bigl\|\sum_{h\in D}u_{n+h}\Bigr\|^{2}
\-\vs{6mm}\\\-
=\frac{1}{|D|^{2}}\frac{1}{|\Phi_{N}|}
\sum_{n\in\Phi_{N}}\sum_{h_{1},h_{2}\in D}\lan u_{n+h_{1}},u_{n+h_{2}}\ran,
}
so
\lequ{
\Bigl(\limsup_{N\ras\infty}\Bigl\|
\frac{1}{|\Phi_{N}|}\sum_{n\in\Phi_{N}}u_{n}\Bigr\|\Bigr)^{2}
=\Bigl(\limsup_{N\ras\infty}\Bigl\|\frac{1}{|D|}\sum_{h\in D}
\frac{1}{|\Phi_{N}|}\sum_{u\in\Phi_{N}}u_{n+h}\Bigr\|\Bigr)^{2}
=\limsup_{N\ras\infty}\Bigl\|\frac{1}{|D|}\sum_{h\in D}
\frac{1}{|\Phi_{N}|}\sum_{n\in\Phi_{N}}u_{n+h}\Bigr\|^{2}
\-\vs{6mm}\\\-
\leq\limsup_{N\ras\infty}\frac{1}{|D|^{2}}\sum_{h_{1},h_{2}\in D}
\frac{1}{|\Phi_{N}|}\sum_{n\in\Phi_{N}}\lan u_{n+h_{1}},u_{n+h_{2}}\ran
\leq\frac{1}{|D|^{2}}\sum_{h_{1},h_{2}\in D}
\limsup_{N\ras\infty}\frac{1}{|\Phi_{N}|}\sum_{n\in\Phi_{N}}\lan u_{n+h_{1}},u_{n+h_{2}}\ran.
}

To get the second assertion, for any finite set $D\sle B$ write
\lequ{
\Bigl|\frac{1}{|D|^{2}}\sum_{h_{1},h_{2}\in D}
\Climsup_{n}\lan u_{n+h_{1}},u_{n+h_{2}}\ran\Bigr|
\-\\\-
\leq\frac{1}{|D|^{2}}\sum_{\sdup{h_{1},h_{2}\in D}{h_{1}\neq h_{2}}}
\bigl|\Climsup_{n}\lan u_{n+h_{1}},u_{n+h_{2}}\ran\bigr|
+\frac{1}{|D|^{2}}\sum_{h\in D}\bigl|\Climsup_{n}\lan u_{n+h},u_{n+h}\ran\bigr|
\leq\eps+\frac{1}{|D|}\sup_{n}\|u_{n}\|^{2}
}
and notice that the second summand tends to zero as $|D|\ras\infty$.%
\endproof
We will also need the following simple ``finitary version'' of the van der Corput trick:
\lemma{P-vdC2}{}{}
Let $u_{1},\ld,u_{N}$ be elements of a Hilbert space.
Then 
$$
\Bigl\|\frac{1}{N}\sum_{n=1}^{N}u_{n}\Bigr\|^{2}
\leq\frac{2}{N}\sum_{h=1}^{N-1}\Bigl|\sum_{n=1}^{N-h}\lan u_{n},u_{n+h}\ran\Bigr|
+\frac{1}{N^{2}}\sum_{n=1}^{N}\|u_{n}\|^{2}.
$$
\endlemma
\proof{}
\lequ{
\Bigl\|\frac{1}{N}\sum_{n=1}^{N}u_{n}\Bigr\|^{2}
=\frac{1}{N^{2}}\sum_{n,m=1}^{N}\lan u_{n},u_{m}\ran
=\frac{1}{N^{2}}\Bigl(\sum_{1\leq n<m\leq N}\lan u_{n},u_{m}\ran
+\sum_{1\leq n<m\leq N}\lan u_{m},u_{n}\ran\Bigr)
+\frac{1}{N^{2}}\sum_{n=1}^{N}\|u_{n}\|^{2}
\-\\\-
=\frac{2}{N^{2}}\Bigl|\sum_{1\leq n<m\leq N}\lan u_{n},u_{m}\ran\Bigr|
+\frac{1}{N^{2}}\sum_{n=1}^{N}\|u_{n}\|^{2}
\leq\frac{2}{N}\sum_{h=1}^{N-1}\Bigl|\sum_{n=1}^{N-h}\lan u_{n},u_{n+h}\ran\Bigr|
+\frac{1}{N^{2}}\sum_{n=1}^{N}\|u_{n}\|^{2}.
}
\frgdsp\endproof
\section{S-APFB}{\BGLFs\ and Besicovitch almost periodicity}

We will now describe and use a ``dynamical'' approach to \BGLFs.
We will focus on functions $\Z\ra\R$.

Let $\cM$ be a torus, $\cM=V/\Gam$,
where $V$ is a finite dimensional $\R$-vector space
and $\Gam$ is a cocompact lattice in $V$,
and let $\pi$ be the projection $V\ra\cM$.
We call {\it a polygon} any bounded subset $P$ of $V$
defined by a system of linear inequalities, strict or non-strict:
$$
P=\Bigl\{v\in V:L_{1}(v)<c_{1},\ld,L_{k}(v)<c_{k},
L_{k+1}(v)\leq c_{k+1},\ld,L_{m}(v)\leq c_{m}\Bigr\},
$$
where $L_{i}$ are linear functions on $V$ and $c_{i}\in\R$.
Let $Q$ be a parallelepiped in $V$ 
such that $\pi\rest{Q}\col Q\ra\cM$ is a bijection.
($Q$ is a fundamental domain of $\cM$ in $V$.)
Assume that $Q=\bigcup_{j=1}^{l}\tP_{j}$ 
is a finite partition of $Q$ into disjoint polygons.
Let a function $\tF$ on $Q$ be the sum, $\tF=L+E$,
of a linear function $L$
and of a function $E$ which is constant on each of $\tP_{j}$.
Finally, let $F$ be the function induced by $\tF$ on $\cM$,
$F=\tF\comp(\pi\rest{Q})^{-1}$.
We will call functions $F$ obtainable this way {\it polygonally broken linear},
or {\it\PGLFs}.

\example{}
The function $\bigl\{2x+\frac{1}{3}\bigr\}$ on $\R/\Z$ is a \PGLF.
\endexample
The following is clear:
\lemma{P-pollin}{}{}
The set of \PGLFs\ on a torus $\cM$
is closed under addition, multiplication by scalars,
and the operation of taking the fractional part.
\endlemma
The following theorem says that \BGLFs\ are dynamically obtainable from \PGLFs:
\theorem{P-dynam}{}{}
For any \BGLF\ $\phi$ there exists a torus $\cM$,
an element $u\in\cM$,
and a \PGLF\ $F$ on $\cM$
such that $\phi(n)=F(nu)$, $n\in\Z$.
\endtheorem
\proof{}
For $\phi(n)=\{an+b\}$, $a,b\in\R$,
take $\cM=\R/\Z$, $u=a\mod\Z$, and $F(x)=\{x+b\}$, $x\in\cM$.

The set of \BGLFs\ satisfying the assertion of the theorem
is closed under addition and multiplication by constants.
Indeed, if a \BGLF\ $\phi$ is represented in the form $\phi(n)=F(nu)$, $n\in\Z$, 
where $F$ is a \PGLF\ on a torus $\cM$ and $u\in\cM$,
then for $a\in\R$ the function $aF$ is a \PGLF\ as well 
and $a\phi(n)=(aF)(nu)$, $n\in\Z$.
If \BGLFs\ $\phi_{1}$, $\phi_{2}$
are represented as $\phi_{1}(n)=F_{1}(nu)$, $\phi_{2}(n)=F_{2}(nv)$, $n\in\Z$,
where $F_{1}$, $F_{2}$ are \PGLFs\ 
on tori $\cM_{1}$ and $\cM_{2}$ respectively,
$u_{1}\in\cM_{1}$ and $u_{2}\in\cM_{2}$,
then the function $F(x_{1},x_{2})=F_{1}(x_{1})+F_{2}(x_{2})$ 
on the torus $\cM_{1}\times\cM_{2}$ is a \PGLF\
and $(\phi_{1}+\phi_{2})(n)=F(n(u_{1},u_{2}))$, $n\in\Z$.

Also, the set of \BGLFs\ satisfying the assertion of the theorem
is closed under the operation of taking the fractional part:
if a \BGLF\ $\phi$ is represented as $\phi(n)=F(nu)$, $n\in\Z$,
where $F$ is a \PGLF\ on a torus $\cM$ and $u\in\cM$,
then the function $\{F\}$ is a \PGLF\
and $\{\phi(n)\}=\{F\}(nu)$, $n\in\Z$.

From the inductive definition of \BGLFs,
it follows that the theorem holds for all \BGLFs.%
\endproof
Any closed subgroup $Z$ of a torus $\cM$ has the form $Z=\cM'\times J$
for some subtorus $\cM'$ of $\cM$ and a finite abelian group $J$.
We will say that a function $F$ on $Z$ is {\it a \PGLF\/}
if the restriction $F\rest{\cM'\times\{i\}}$ 
is a \PGLF\ on the torus $\cM'\times\{i\}$ for every $i\in J$.

If $Z$ is a closed subgroup of a torus $\cM$
and $F$ is a \PGLF\ on $\cM$,
then $F\rest{Z}$ is a \PGLF\ on $Z$.
In the environment of \rfr{P-dynam},
putting $Z=\ovr{\Z u}$,
we obtain the following:
\proposition{P-edynam}{}{}
For any \BGLF\ $\phi$ there exists a compact abelian group $Z$,
of the form $Z=\cM'\times J$, where $\cM'$ is a torus and $J$ is a finite cyclic group,
an element $u\in Z$, whose orbit $\Z u$ is dense 
(and so, uniformly distributed)
in $Z$,
and a \PGLF\ $F$ on $Z$
such that $\phi(n)=F(nu)$, $n\in\Z$.
\endproposition
\corollary{P-dens}{}{}
For any \BGLF\ $\phi$, the limit $\Clim_{n}\phi(n)$ exists.
For any \BGLFs\ $\phi_{1},\ld,\phi_{k}$, for $\phi=(\phi_{1},\ld,\phi_{k})$,
and for any polygon $P\sle\R^{k}$,
the density of the set $\{n\in\Z:\phi(n)\in P\}$ exists.
\endcorollary
As another corollary of \rfr{P-edynam}, we get the following result:
\proposition{P-Besic}{}{}
Let $\phi\col\Z\ra\R$ be a \BGLF.
For any $\eps>0$ 
there exists $h\in\Z$ 
such that $\Dlimsup_{n}|\phi(n+h)-\phi(n)|<\eps$,
and there exists a trigonometric polynomial $q$
such that $\Dlimsup_{n}|\phi(n)-q(n)|<\eps$.
\endproposition
\remark{}
Functions with these properties are called {\it Besicovitch almost periodic}
(at least, in the case the F{\o}lner sequence 
with respect to which the densities are measured 
is $\Phi_{N}=[-N,N]$, $N\in\N$).
Any function obtainable dynamically 
with the help of a rotation of a compact commutative Lie group 
and a Riemann integrable function thereon
is such.
\endremark
\proof{}
Represent $\phi$ in the form $\phi(n)=F(nu)$, $n\in\Z$, as in \rfr{P-edynam}.
Let $Z=\bigcup_{j=1}^{l}P_{j}$ be the polygonal partition of $Z$
such that $F$ is linear on each of $P_{j}$.
Let $U$ be a $\del$-neighborhood of $\bigcup_{j=1}^{l}\partial P_{j}$
with $\del>0$ small enough so that $\lam(U)<\eps$,
where $\lam$ is the normalized Haar measure on $Z$.
Let $\hF$ be a continuous function on $Z$ which coincides with $F$ on $Z\sm U$
and such that $\sup|\hF|\leq\sup|F|=\sup|\phi|$.
Let $\hphi(n)=\hF(nu)$, $n\in\Z$.
The sequence $(nu)$ is uniformly distributed on $Z$,
thus $\Du\bigl(\{n\in\Z:nu\in U\}\bigr)<\eps$,
and so $\Du\bigl(\{n:\phi(n)\neq\hphi(n)\}\bigr)=\Du\bigl(\{n:F(nu)\neq\hF(nu)\}\bigr)<\eps$.
Since $\hF$ is uniformly continuous,
for any $h\in\Z$ for which $hu$ is close enough to 0
we have $\bigl|\hF(v+hu)-\hF(v)\bigr|<\eps$ for all $v\in Z$,
so $|\hphi(n+h)-\hphi(n)|<\eps$ for all $n\in\Z$.
This implies that $\Dlimsup_{n}|\phi(n+h)-\phi(n)|<\eps+2\eps\sup|\phi|$.
And if $\Theta$ is a finite linear combination of characters of $Z$ 
such that $|\hF-\Theta|<\eps$,
then for the trigonometric polynomial $q(n)=\Theta(nu)$, $n\in\Z$,
we have $|\hphi(n)-q(n)|<\eps$ for all $n$,
which implies that $\Dlimsup_{n}|\phi(n)-q(n)|<\eps+(\sup|\phi|+\sup|q|)\eps
=\eps+(2\sup|\phi|+\eps)\eps$.%
\endproof
\corollary{P-APset}{}{}
If $\phi$ is a \BGLF\ $\Z\ra\Z$,
then for any $\eps>0$ there exists $h\in\Z$ 
such that $\D\bigl(\bigl\{n\in\Z:\phi(n+h)=\phi(n)\bigr\}\bigr)>1-\eps$.
\endcorollary
\npar
(Notice that the density of the set $\bigl\{n\in\Z:\phi(n+h)=\phi(n)\bigr\}$
exists by \rfr{P-dens}.)
We now turn to unbounded \GLFs.
From \rfr{P-LBP} and \rfr{P-dynam} 
we see that any \GLF\ $\phi$ is representable in the form $\phi(n)=an+F(nu)$,
where $a\in\R$, $F$ is a \PGLF\ on a torus $\cM$, and $u\in\cM$.
Given several \GLFs\ $\phi_{1},\ld,\phi_{k}$, we can read them off a single torus:
for each $i$ represent $\phi_{i}$ in the form $\phi_{i}(n)=a_{i}n+F_{i}(nu_{i})$,
where $a_{i}\in\R$, 
$F_{i}$ is a \PGLF\ on a torus $\cM_{i}$,
and $u_{i}\in\cM_{i}$,
put $\cM=\prod_{i=1}^{k}\cM_{i}$, $u=(u_{1},\ld,u_{k})\in\cM$,
and lift $F_{1},\ld,F_{k}$ to a function on $\cM$;
then $\phi_{i}(n)=a_{i}n+F_{i}(nu)$, $n\in\Z$, $i=1,\ld,k$.
As a corollay, we get:
\proposition{P-Part}{}{}
Given \GLFs\ $\phi_{1},\ld,\phi_{k}$,
there exists a torus $\cM$, an element $u\in\cM$, 
and a polygonal partition $\cM=\bigcup_{j=1}^{l}P_{j}$,
such that for each $i$, $j$, $\phi_{i}(n+h)-\phi_{i}(n)$ does not depend on $n$ 
if both $nu,(n+h)u\in P_{j}$.
\endproposition
\proof{}
Let $\cM$, $u$, and $F_{i}$ be as above;
let $\cM=V/\Gam$ where $V$ is a vector space and $\Gam$ is a lattice in $V$,
$\pi$ be the projection $V\ra\cM$, 
$Q\sln V$ be the fundamental domain of $\cM$ in $V$,
and $\tF_{i}=F\comp\pi\rest{Q}$, $i=1,\ld,k$.
Choose a partition $\cM=\bigcup_{j=1}^{l}P_{j}$ of $\cM$
such that for each $j$ and each $i$,
the function $F_{i}$ is linear on $P_{j}$,
and, additionally, for each $j$, 
$\bigl((\tP_{j}-\tP_{j})-(\tP_{j}-\tP_{j})\bigr)\cap\Gam=\{0\}$,
where $\tP_{j}=\pi^{-1}(P_{j})\cap Q$.
Then for any $i$ and $j$,
for $v,w\in P_{j}$, $F_{i}(v)-F_{i}(w)$ depends on $v-w$ only.
Indeed, let $v_{1},w_{1},v_{2},w_{2}\in P_{j}$ be such that $v_{1}-w_{1}=v_{2}-w_{2}$;
let $\tv_{t}=\pi\rest{Q}^{-1}(v_{t})$, $\tw_{t}=\pi\rest{Q}^{-1}(w_{t})$, $t=1,2$,
then $(\tv_{1}-\tw_{1})-(\tv_{2}-\tw_{2})\in\Gam$, so $=0$,
thus
$$
F_{i}(v_{1})-F_{i}(w_{1})=\tF_{i}(\tv_{1})-\tF_{i}(\tw_{1})=L_{i}(\tv_{1})-L_{i}(\tw_{1})
=L_{i}(\tv_{1}-\tw_{1})=L_{i}(\tv_{2}-\tw_{2})=F_{i}(v_{2})-F_{i}(w_{2}),
$$
where $L_{i}$ is the linear function on $V$ 
that coincides with $\tF_{i}$ on $\tP_{j}$ up to a constant.
Now, if $n$ and $h$ are such that both $nu,(n+h)u\in P_{j}$ for some $j$,
then for any $i$,
$\phi_{i}(n+h)-\phi_{i}(n)=a_{i}h+F_{i}(nu+hu)-F_{i}(nu)$,
and $F_{i}(nu+hu)-F_{i}(nu)$ does not depend on $n$.%
\endproof
A set $H\sle\Z$ is said to be {\it a Bohr set\/}
if $H$ contains a nonempty subset of the form $\{n\in\Z:nu\in W\}$,
where $u$ and $W$ are an element and an open subset of a torus.
Any Bohr set is infinite and has positive density 
(with respect to any F{\o}lner sequence in $\Z$).
The following proposition says that (several) \GLFs\ are ``almost linear'' 
along a Bohr set:
\proposition{P-Hlin}{}{}
For any \GLFs\ $\phi_{1},\ld,\phi_{k}$ and any $\eps>0$ 
there exists a Bohr set $H\sle\Z$ and constants $C_{1},\ld,C_{k}$
such that for any $h\in H$,
$$
\D\Bigl(\Bigl\{n\in\Z:\phi_{i}(n+h)=\phi_{i}(n)+\phi_{i}(h)+C_{i},\ 
i=1,\ld,k\Bigr\}\Bigr)>1-\eps.
$$
\endproposition
\proof{}
First of all, for any $h$, the density of the set 
$\bigl\{n\in\Z:\phi_{i}(n+h)=\phi_{i}(n)+\phi_{i}(h)+C_{i},\ i=1,\ld,k\bigr\}$
exists by \rfr{P-dens}.

Let $\cM$ be a torus, $u\in\cM$, 
and $F_{1},\ld,F_{k}$ be \PGLFs\ on $\cM$
such that $\phi_{i}(n)=F_{i}(nu)$, $i=1,\ld,k$.
Let $Z=\ovr{\Z u}$;
then the sequence $(nu)_{n\in\Z}$ is uniformly distributed in $Z$.
Let $Z=\bigcup_{j=1}^{l}P_{j}$ be the polygonal partition of $Z$
such that for every $i$ and $j$, $F_{i}\rest{P_{j}}=L_{i}+C_{i,j}$,
where $L_{i}$ is linear and $C_{i,j}$ is a constant.
Let $\del>0$ be small enough 
so that $\lam(U)<\eps$
where $U$ is the $\del$-neighborhood of the set $\bigcup_{j=1}^{l}\partial P_{j}$
and $\lam$ is the normalized Haar measure on $Z$.
Let $W_{0}$ be the $\del$-neighborhood of 0.
Now, for any $w\in W_{0}$,
$\D\bigl(\bigl\{n\in\Z:
nu\in P_{j_{1}},\ nu+w\in P_{j_{2}},\ j_{1}\neq j_{2}\bigr\}\bigr)<\eps$. 
Choose $j_{0}$ for which 0 is a limit point 
of the interior $\Po_{j_{0}}$ of $P_{j_{0}}$,
let $W=\Po_{j_{0}}\cap W_{0}$, 
and let $H=\{n\in\Z:nu\in W\}$.
Then for any $w\in W$ and any $i$,
whenever $v,v+w\in P_{j}$ for some $j$ we have
$$
F_{i}(v+w)=L_{i}(v+w)+C_{i,j}=L_{i}(v)+C_{i,j}+L_{i}(w)+C_{i,j_{0}}-C_{i,j_{0}}
=F_{i}(v)+F_{i}(w)+C_{i},
$$
where $C_{i}=-C_{i,j_{0}}$.
For any $h\in H$ let $E_{h}=\bigl\{n\in\Z:\hbox{$nu,(n+h)u\in P_{j}$ for some $j$}\bigr\}$;
then $\D(E_{h})>1-\eps$,
and for any $n\in E_{h}$ and any $i$, $\phi_{i}(n+h)=\phi_{i}(n)+\phi_{i}(h)+C_{i}$.%
\endproof
\section{S-GLST}{Generalized linear sequences of transformations}

{\it A generalized linear sequence (a \GLS)\/} in a commutative group $G$ 
is a sequence of the form 
$\T(n)=T_{1}^{\phi_{1}(n)}\cds T_{r}^{\phi_{r}(n)}$, $n\in\Z$,
where $T_{1},\ld,T_{r}\in G$ and $\phi_{1},\ld,\phi_{r}$ are \GLFs\ $\Z\ra\Z$.
We say that $\T$ is a \BGLS\ if $\phi_{1},\ld,\phi_{r}$ are \BGLFs.
\rfr{P-Intp}, \rfr{P-APset}, \rfr{P-Part}, and \rfr{P-Hlin}
imply the following properties of \GLSs:
\proposition{P-perT}{}{}
Let $G$ be a commutative group.
\npar{\rm(i)}
If $\T$ is a \GLS\ in $G$,
then for any $h\in\Z$ the sequence $\T(n)^{-1}\T(n+h)$, $n\in\Z$, is a \BGLS.
\npar{\rm(ii)}
If $\T$ is a \BGLS\ in $G$,
then for any $\eps>0$ there exists $h\in\Z$ 
such that $\Dl\bigl(\bigl\{n\in\Z:\T(n+h)=\T(n)\bigr\}\bigr)>1-\eps$.
\npar{\rm(iii)}
If $\T_{1},\ld,\T_{k}$ are \GLSs\ in $G$
(or, more generally, in distinct commutative groups $G_{1},\ld,G_{k}$ respectively)
then there exist a torus $\cM$, an element $u\in\cM$,
and a polygonal partition $\cM=\bigcup_{j=1}^{l}P_{j}$
such that for any $i$, $j$,
$\T_{i}(n)^{-1}\T_{i}(n+h)$ does not depend on $n$ whenever $nu,(n+h)u\in P_{j}$.
\npar{\rm(iv)}
If $\T_{1},\ld,\T_{k}$ are \GLSs\ in $G$
(or, more generally, in distinct commutative groups $G_{1},\ld,G_{k}$ respectively)
then for any $\eps>0$
there exist a Bohr set $H\sle\Z$ and elements $S_{1},\ld,S_{k}\in G$
(respectively, $S_{i}\in G_{i}$, $i=1,\ld,k$)
such that for any $h\in H$ 
the set $E_{h}=\bigl\{n\in\Z:\T_{i}(n+h)=\T_{i}(n)\T_{i}(h)S_{i},\ i=1,\ld,k\bigr\}$ 
satisfies $\Dl(E_{h})>1-\eps$.
\endproposition
If $\T$ is a \GLS\ of unitary operators on a Hilbert space $\H$,
then via the spectral theorem, \rfr{P-dens} implies the following:
\lemma{P-clim}{}{}
For any $f\in\H$, $\Clim_{n}\T(n)f$ exists.
\endlemma
We now fix a commutative group $G$ of measure preserving transformations 
of a probability measure space $(X,\mu)$,
and denote by $\gT$ the set of \GLSs\ of transformations in $G$.

\ndefinition{D-erg}{}
If $\T$ is a sequence of measure preserving transformations of $X$
(or just a sequence of unitary operators on a Hilbert space $\H$),
we say that $\T$ is {\it ergodic\/} 
if $\Clim_{n}\T(n)f=\int_{X}f\,d\mu$ for all $f\in L^{2}(X)$
(respectively, $\Clim_{n}\T(n)f=0$ for all $f\in\H$).
We will also say that $\T$ is {\it weakly mixing\/}
if for any $f,g\in L^{2}(X)$ 
one has $\Dlim_{n}\int_{X}\T(n)f\cd g\,d\mu=\int_{X}f\,d\mu\int_{X}g\,d\mu$
(respectively, $\Dlim_{n}\lan\T(n)f,g\ran=0$ for all $f,g\in\H$).
\endndefinition
\nremark{R-indep}{}
We have defined our $\Clim${\kern1pt}s, and so, ergodicity of a sequence of transformations,
with respect to a fixed F{\o}lner sequence in $\Z$.
However, since,
for any \GLS\ $\T$ of measure preserving transformations, or of unitary operators, 
and for any (function or vector) $f$,
$\Clim\T(n)f$ exists with respect to any F{\o}lner sequence,
this limit is the same for all F{\o}lner sequences;
thus, the ergodicity of \GLSs\ does not depend on the choice of the F{\o}lner sequence.
\endnremark

Let $\Hc\oplus\Hwm$ be the compact/weak~mixing decomposition of $L^{2}(X)$ induced by $G$,
meaning that $\Hc$ is the subspace of $L^{2}(X)$ 
on which all elements of $G$ act in a compact way
and $\Hwm$ is the orthocomplement of $\Hc$;
then for any $g\in\Hwm$ there exists a transformation $T\in G$
that acts on $g$ in a weakly mixing fashion.
Notice also that if $T\in G$ is ergodic, 
then $T$ is weakly mixing on $\Hwm$.
The following theorem says
that any ergodic sequence from $\gT$ is weakly mixing on $\Hwm$:
\theorem{P-cwm}{}{}
If $\T\in\gT$ is ergodic,
then for any $f\in\Hwm$ and $g\in L^{2}(X)$ 
one has $\Dlim_{n}\int_{X}\T(n)f\cd g\,d\mu=0$.
\endtheorem
We first prove that $\T$ has no ``eigenfunctions'' in $\Hwm$:
\lemma{P-noeig}{}{}
If $\T\in\gT$ is ergodic,
then for any $f\in\Hwm$ and $\lam\in\C$ with $|\lam|=1$
one has $\Clim_{n}\lam^{n}\T(n)f=0$.
\endlemma
\proof{}
We may and will assume that $|f|\leq1$.
Fix $g\in\Hwm$ with $|g|\leq 1$, 
and let $T\in G$ be a transformation that acts weakly mixingly on $g$.
We are going to apply the van der Corput trick (\rfr{P-vdC} above)
to the sequence $f_{n}=\lam^{n}T^{n}\T(n)f\cd T^{n}g$, $n\in\Z$.
Let $\eps>0$,
and let a Bohr set $H\sle\Z$, a transformation $S\in G$, and sets $E_{h}\sle\Z$, $h\in H$,
be as in \rfr{P-perT}(iv), applied to the single \GLS\ $\T$.
Let $h_{1},h_{2}\in H$;
for any $n\in E_{h_{1}}\cap E_{h_{2}}$ one has
\lequ{
\lan f_{n+h_{1}},f_{n+h_{2}}\ran
=\int_{X}f_{n+h_{1}}\af_{n+h_{2}}\,d\mu
\-\\\-
=\int_{X}\lam^{h_{1}}T^{n+h_{1}}\T(n+h_{1})f\cd T^{n+h_{1}}g\cd 
\bar{\lam}^{h_{2}}T^{n+h_{2}}\T(n+h_{2})\af\cd T^{n+h_{2}}\ag\,d\mu
\-\\\-
=\int_{X}\T(n)\bigl(\lam^{h_{1}-h_{2}}T^{h_{1}}\T(h_{1})Sf\cd T^{h_{2}}\T(h_{2})S\af\bigr)
\cd(T^{h_{1}}g\cd T^{h_{2}}\ag)\,d\mu
\\\-
=\int_{X}\T(n)\tf_{h_{1},h_{2}}\cd(T^{h_{1}}g\cd T^{h_{2}}\ag)\,d\mu,
}
where $\tf_{h_{1},h_{2}}=\lam^{h_{1}-h_{2}}T^{h_{1}}\T(h_{1})Sf\cd T^{h_{2}}\T(h_{2})S\af$.
Since $\T$ is ergodic,
\equ{
\Clim_{n}\int_{X}\T(n)\tf_{h_{1},h_{2}}\cd(T^{h_{1}}g\cd T^{h_{2}}\ag)\,d\mu
=\Bigl(\int_{X}\tf_{h_{1},h_{2}}d\mu\Bigr)
\Bigl(\int_{X}T^{h_{1}}g\cd T^{h_{2}}\ag\,d\mu\Bigr),
}
and since $\Dl(E_{h_{1}}\cap E_{h_{2}})>1-2\eps$ and $|\tf_{h_{1},h_{2}}|,|g|\leq1$,
$$
\Bigl|\Climsup_{n}\lan f_{n+h_{1}},f_{n+h_{2}}\ran\Bigr|
\leq\Bigl|\int_{X}T^{h_{1}}g\cd T^{h_{2}}\ag\,d\mu\Bigr|+2\eps.
$$
Since $\Dlim_{h}\int_{X}T^{h}g\cd g'\,d\mu=0$ for any $g'\in L^{2}(X)$, 
and since $\Dl(H)>0$,
we can construct an infinite set $B\sle H$ such that
$\bigl|\int_{X}T^{h_{1}}g\cd T^{h_{2}}\ag\,d\mu\bigr|<\eps$ 
for any distinct $h_{1},h_{2}\in B$.
Then for any distinct $h_{1},h_{2}\in B$ we have
$\bigl|\Climsup_{n}\lan f_{n+h_{1}},f_{n+h_{2}}\ran\bigr|<3\eps$,
which, by \rfr{P-vdC}, implies that $\Climsup_{\|\cd\|,n}f_{n}\leq\sqrt{3\eps}$.
Since $\eps$ is arbitrary, $\Clim_{n}f_{n}=0$.

Now, let $\hf=\Clim_{n}\lam^{n}\T(n)f\in\Hwm$.
Then for any $g\in L^{2}(X)$,
$$
\int_{X}\hf\cd g\,d\mu=\Clim_{n}\int_{X}\lam^{n}\T(n)f\cd g\,d\mu=\Clim_{n}\int_{X}f_{n}d\mu=0.
$$
Hence, $\hf=0$.%
\endproof
\proof{of \rfr{P-cwm}}
Let $\T(n)=T_{1}^{\phi_{1}(n)}\cds T_{r}^{\phi_{r}(n)}$, $n\in\Z$,
where $T_{i}\in G$ and $\phi_{i}$ are \GLFs.
By \rfr{P-LBP}, for each $j=1,\ld,r$
one has $\phi_{j}(n)=a_{j}n+\psi_{j}(n)$, $n\in\Z$,
where $a_{j}\in\R$ and $\psi_{j}$ is a \BGLF.
Considering $T_{1},\ld,T_{r}$ as unitary operators on $\Hwm$,
immerse them into commuting continuous unitary flows $(T_{i}^{t})_{t\in\R}$,
and let $T=T_{1}^{a_{1}}\cds T_{r}^{a_{r}}$.

Based on \rfr{P-noeig},
we are going to show that $T$ has no eigenvectors in $\Hwm$.
Assume, in the course of contradiction,
that there exists $f\in\Hwm$, with $\|f\|=1$, such that $Tf=\lam f$.
Let $\S(n)=T_{1}^{\psi_{1}(n)}\cds T_{r}^{\psi_{r}(n)}$,
so that $\T(n)=T^{n}\S(n)$, $n\in\Z$.
Fix $\eps>0$.
Let $I$ be an interval in $\R^{r}$ 
that contains the range $\psi(\Z)$ of the function $\psi=(\psi_{1},\ld,\psi_{r})$.
Partition $I$ to subintervals $I_{1},\ld,I_{l}$ small enough
so that for each $i=1,\ld,l$, 
for some $f_{i}\in\Hwm$ 
one has $(T_{1}^{z_{1}}\cds T_{r}^{z_{r}})^{-1}f\app^{\eps}f_{i}$ 
for all $(z_{1},\ld,z_{r})\in I_{i}$.
(Here and below, for $g_{1},g_{2}\in L^{2}(X)$,
``$g_{1}\app^{\eps}g_{2}$'' means that $\|g_{1}-g_{2}\|\leq\eps$.)
For each $i=1,\ld,l$ let $A_{i}=\bigl\{n:\psi(n)\in I_{i}\bigr\}$;
by \rfr{P-indI}, the indicator function $1_{A_{i}}$ is a \UGLF,
and thus by \rfr{P-Besic} there exists a trigonometric polynomial $q_{i}$
such that $\Dlimsup_{n}|1_{A_{i}}(n)-q_{i}(n)|<\eps/l$.
For any $i$, for any $n\in A_{i}$, 
$T^{n}f=\T(n)\S(n)^{-1}f\app^{\eps}\T(n)f_{i}$,
thus,
\equ{
\D(A_{i})f=\Clim_{n}1_{A_{i}}(n)\lam^{-n}T^{n}f
\app^{\sD(A_{i})\eps}\Clim_{n}1_{A_{i}}(n)\lam^{-n}\T(n)f_{i}
\app^{\eps/l}\Clim_{n}q_{i}(n)\lam^{-n}\T(n)f_{i}.
}
By \rfr{P-noeig}, the last limit is equal to 0;
summing this up for $i=1,\ld,l$ we get $f\app^{2\eps}0$.
Since $\eps$ is arbitrary, $f=0$.

Hence, $T$ is weakly mixing on $\Hwm$,
that is, for any $f,g\in\Hwm$, $\Dlim_{n}\lan T^{n}f,g\ran=0$.
Let $f,g\in\Hwm$ and let $\eps>0$.
The set $\{\S(n)^{-1}g,\ n\in\Z\}$ is totally bounded;
let $\{g_{1},\ld,g_{k}\}$ be an $\eps$-net in this set.
Then 
$$
\Dlimsup_{n}\bigl|\lan\T(n)f,g\ran\bigr|
=\Dlimsup_{n}\bigl|\lan T^{n}f,\S(n)^{-1}g\ran\bigr|
<\Dlim_{n}\max_{i}|\lan T^{n}f,g_{i}\ran|+\eps=\eps.
$$
Since $\eps$ is arbitrary, $\Dlim_{n}\lan\T(n)f,g\ran=0$.%
\endproof
\section{S-JE}{Joint ergodicity of several \GLSs\ of transformations}

We now start dealing with several \GLSs\ of measure preserving transformations.
We preserve the notations $G$, $\gT$, $\Hc$, and $\Hwm$ from the preceding section.

Given functions $f_{1},\ld,f_{k}$ on $X$,
the tensor product $\bigotimes_{i=1}^{k}f_{i}=f_{1}\tens\cds\tens f_{k}$ 
is the function $f(x_{1},\ld,x_{k})=f_{1}(x_{1})\cds f_{k}(x_{k})$ on $X^{k}$
(whereas the product $\prod_{i=1}^{k}f_{i}=f_{1}\cds f_{k}$ 
is the function $f_{1}(x)\cds f_{k}(x)$ on $X$).
\lemma{P-tens}{}{}
If $\T_{1},\ld,\T_{k}\in\gT$ are ergodic,
then for any functions $f_{1},\ld,f_{k}\in L^{\infty}(X)$ 
with $f_{i}\in\Hwm$ for at least one $i$,
$\Clim_{n}\bigotimes_{i=1}^{k}\T_{i}(n)f_{i}=0$ in $L^{2}(X^{k})$.
\endlemma
\proof{}
Assume that $f_{1}\in\Hwm$.
Let $\hf=\Clim_{n}\bigotimes_{i=1}^{k}\T_{i}(n)f_{i}$.
For any $g_{1},\ld,g_{k}\in L^{\infty}(X)$ we have
\equ{
\Blan\hf,\bigotimes_{i=1}^{k}g_{i}\Bran
=\Clim_{n}\int_{X^{k}}\bigotimes_{i=1}^{k}\T_{i}(n)f_{i}\cd\bigotimes_{i=1}^{k}\ag_{i}\,d\mu^{k}
=\Clim_{n}\prod_{i=1}^{k}\int_{X}\T_{i}(n)f_{i}\cd\ag_{i}\,d\mu
=0
}
since $\Dlim_{n}\int_{X}\T_{1}(n)f_{1}\cd\ag_{1}\,d\mu=0$ by \rfr{P-cwm}.
Hence, $\hf=0$.%
\endproof
Given transformations $T_{1},\ld,T_{k}$ of $X$,
$T_{1}\times\cds\times T_{k}$ is the transformation of $X^{k}$
defined by $(T_{1}\times\cds\times T_{k})(x_{1},\ld,x_{k})
=\bigl(T_{1}x_{1},\ld,T_{k}x_{k}\bigr)$.
Notice that if $\T_{1},\ld,\T_{k}$ are sequences of transformations of $X$
such that $\T_{1}\times\cds\times\T_{k}$ is ergodic,
then $\T_{1},\ld,\T_{k}$ are ergodic,
and, moreover, $\T_{i_{1}}\times\cds\times\T_{i_{l}}$ is ergodic 
for any $1\leq i_{1}<\ld<i_{l}\leq k$.
\lemma{P-induc}{}{}
Let $\T_{1},\ld,\T_{k}\in\gT$ be such that
the \GLSs\ $\T_{1}^{-1}\T_{2},\ld,\T_{1}^{-1}\T_{k}$ of transformations of $X$ are ergodic
and the \GLS\ $\T_{1}\times\cds\times\T_{k}$ of transformations of $X^{k}$ is ergodic.
Then the \GLS\ $(\T_{1}^{-1}\T_{2})\times\cds\times(\T_{1}^{-1}\T_{k})$ 
of transformations of $X^{k-1}$ is also ergodic.
\endlemma
\proof{}
Since the span of the functions of the form $f_{2}\otimes\cds\otimes f_{k}$ 
with $f_{2},\ld,f_{k}\in L^{\infty}(X)$ is dense in $L^{2}(X^{k-1})$,
it suffices to show that for any $f_{2},\ld,f_{k}\in L^{\infty}(X)$
with $\int_{X}f_{2}\,d\mu=\ld=\int_{X}f_{k}\,d\mu=0$ one has
\equn{f-induc}{
\Clim_{n}\bigotimes_{i=2}^{k}\T_{1}(n)^{-1}\T_{i}(n)f_{i}=0
}
in $L^{2}(X^{k-1})$.
If at least one of $f_{i}$ is in $\Hwm$, this is true by \rfr{P-tens}.
If $f_{2},\ld,f_{k}\in\Hc$,
we may assume that $f_{2},\ld,f_{k}$ are nonconstant eigenfunctions of the elements of $G$,
so that $\T_{1}(n)f_{i}=\tau_{i}(n)f_{i}$ and $\T_{i}(n)f_{i}=\lam_{i}(n)f_{i}$, $i=2,\ld,k$,
for some (multiplicative) \GLSs\ $\tau_{i}$, $\lam_{i}$ in $\{z\in\C:|z|=1\}$.
Put $f_{1}=\ovr{f_{2}\cds f_{k}}$.
Since $\T_{1}\times\cds\times\T_{k}$ is ergodic,
we have $\Clim_{n}\bigotimes_{i=1}^{k}\T_{i}(n)f_{i}=0$,
which implies that $\Clim_{n}\prod_{i=1}^{k}\ovr{\tau_{i}(n)}\lam_{i}(n)=0$,
which then implies \frfr{f-induc}.%
\endproof
\ndefinition{D-jerg}{}
We say that sequences $\T_{1},\ld,\T_{k}$ of measure preserving transformations of $X$
are {\it jointly ergodic\/}
if for any $f_{1},\ld,f_{k}\in L^{\infty}(X)$,
$\Clim_{n}\prod_{i=1}^{k}\T_{i}(n)f_{i}=\prod_{i=1}^{k}\int_{X}f_{i}\,d\mu$ in $L^{2}(X)$.
\endndefinition
Notice that if $\T_{1},\ld,\T_{k}$ are jointly egodic,
then $\T_{1},\ld,\T_{k}$ are ergodic,
and, moreover, $\T_{i_{1}},\ld,\T_{i_{l}}$ are jointly ergodic
for any $1\leq i_{1}<\ld<i_{l}\leq k$.

We are now in position to prove our main result:
\theorem{P-joint}{}{}
\GLSs\ $\T_{1},\ld,\T_{k}\in\gT$ are jointly ergodic
iff the \GLSs\ $\T_{i}^{-1}\T_{j}$ are ergodic for all $i\neq j$
and the \GLS\ $\T_{1}\times\cds\times\T_{k}$ is ergodic.
\endtheorem
\proof{}
Assume that $\T_{1},\ld,\T_{k}\in\gT$ are jointly ergodic.
Let $i,j\in\{1,\ld,k\}$, $i\neq j$, let $f\in L^{\infty}(X)$,
and let $\hf=\Clim_{n}\T_{i}^{-1}(n)\T_{j}(n)f$.
Then for any $g\in L^{\infty}(X)$ we have
$$
\lan\hf,g\ran
=\Clim_{n}\int_{X}\T_{i}^{-1}(n)\T_{j}(n)f\cd\ag\,d\mu
=\int_{X}\Clim_{n}\T_{j}(n)f\cd\T_{i}(n)\ag\,d\mu
=\int_{X}f\,d\mu\int_{X}\ag\,d\mu.
$$
Hence, $\hf=\int_{X}f\,d\mu$, so $\T_{i}^{-1}\T_{j}$ is ergodic.

To prove that $\T_{1}\times\cds\times\T_{k}$ is also ergodic,
it suffices to show that for any $f_{1},\ld,f_{k}\in L^{\infty}(X)$
with $\int_{X}f_{1}\,d\mu=\ld=\int_{X}f_{k}\,d\mu=0$,
$\Clim_{n}\bigotimes_{i=1}^{k}\T_{i}(n)f_{i}=0$.
If at least one of $f_{i}$ is in $\Hwm$,
this is true by \rfr{P-tens}.
If $f_{i}\in\Hc$ for all $i$,
we may assume that $f_{1},\ld,f_{k}$ are nonconstant eigenfunctions of the elements of $G$,
and then $\T_{i}(n)f_{i}=\lam_{i}(n)f_{i}$, $i=1,\ld,k$,
for some \GLSs\ $\lam_{1},\ld,\lam_{k}$ in $\{z\in\C:|z|=1\}$.
In this case both 
$\bigotimes_{i=1}^{k}\T_{i}(n)f_{i}=\lam(n)\bigotimes_{i=1}^{k}f_{i}$
and $\prod_{i=1}^{k}\T_{i}(n)f_{i}=\lam(n)\prod_{i=1}^{k}f_{i}$,
where $\lam(n)=\prod_{i=1}^{k}\lam_{i}(n)$, $n\in\Z$.
Since $\Clim_{n}\prod_{i=1}^{k}\T_{i}(n)f_{i}=0$, 
we have $\Clim_{n}\lam(n)=0$,
and so, $\Clim_{n}\bigotimes_{i=1}^{k}\T_{i}(n)f_{i}=0$.

\vbreak{-9000}{2mm}
Conversely, assume that $\T_{i}^{-1}\T_{j}$ are ergodic for all $i\neq j$
and $\T_{1}\times\cds\times\T_{k}$ is ergodic,
and let $f_{1},\ld,f_{k}\in L^{\infty}(X)$.
If all $f_{i}\in\Hc$, then, again,
we may assume that $f_{1},\ld,f_{k}$ are nonconstant eigenfunctions of the elements of $G$,
so $\bigotimes_{i=1}^{k}\T_{i}(n)f_{i}=\lam(n)\bigotimes_{i=1}^{k}f_{i}$
and $\prod_{i=1}^{k}\T_{i}(n)f_{i}=\lam(n)\prod_{i=1}^{k}f_{i}$,
and since now $\Clim_{n}\bigotimes_{i=1}^{k}\T_{i}(n)f_{i}=0$,
we obtain that $\Clim_{n}\prod_{i=1}^{k}\T_{i}(n)f_{i}=0$ as well.

It remains to show that $\Clim_{n}\prod_{i=1}^{k}\T_{i}(n)f_{i}=0$
whenever $f_{i}\in\Hwm$ for at least one $i$.
We will assume that $f_{1}\in\Hwm$ and that $|f_{i}|\leq 1$ for all $i$.
We will use the van der Corput trick and induction on $k$.
Let $\eps>0$.
Let a Bohr set $H\sle\Z$, transformations $S_{1},\ld,S_{k}\in G$,
and sets $E_{h}\sle\Z$, $h\in H$, be as in \rfr{P-perT}(iv),
applied to the \GLSs\ $\T_{1},\ld,\T_{k}$.
Let $h_{1},h_{2}\in H$;
for any $n\in E_{h_{1}}\cap E_{h_{2}}$ we have
\lequ{
\Blan\prod_{i=1}^{k}\T_{i}(n+h_{1})f_{i},\prod_{i=1}^{k}\T_{i}(n+h_{2})f_{i}\Bran
=\int_{X}\prod_{i=1}^{k}\T_{i}(n+h_{1})f_{i}\cd\prod_{i=1}^{k}\T_{i}(n+h_{2})\af_{i}\,d\mu
\-\\\-
=\int_{X}\bigl(\T_{1}(h_{1})S_{1}f_{1}\cd\T_{1}(h_{2})S_{1}\af_{1}\bigr)\cd
\prod_{i=2}^{k}\T_{1}^{-1}(n)\T_{i}(n)
\bigl(\T_{i}(h_{1})S_{i}f_{i}\cd\T_{i}(h_{2})S_{i}\af_{i}\bigr)d\mu.
}
By \rfr{P-induc}, $\T_{1}^{-1}\T_{2}\times\cds\times\T_{1}^{-1}\T_{k}$ is ergodic,
thus by induction on $k$,
$\T_{1}^{-1}\T_{2},\ld,\T_{1}^{-1}\T_{k}$ are jointly ergodic,
so
\lequ{
\Clim_{n}\int_{X}\bigl(\T_{1}(h_{1})S_{1}f_{1}\cd\T_{1}(h_{2})S_{1}\af_{1}\bigr)\cd
\prod_{i=2}^{k}\T_{1}^{-1}(n)\T_{i}(n)
\bigl(\T_{i}(h_{1})S_{i}f_{i}\cd\T_{i}(h_{2})S_{i}\af_{i}\bigr)d\mu
\-\\\-
=\prod_{i=1}^{k}\int_{X}\T_{i}(h_{1})S_{i}f_{i}\cd\T_{i}(h_{2})S_{i}\af_{i}\,d\mu.
}
Since $\Dl(E_{h_{1}}\cap E_{h_{2}})>1-2\eps$ and $|f_{1}|,\ld,|f_{k}|\leq 1$, 
we get
\equ{
\Bigl|\Climsup_{n}\Blan\prod_{i=1}^{k}\T_{i}(n+h_{1})f_{i},
\prod_{i=1}^{k}\T_{i}(n+h_{2})f_{i}\Bran\Bigr|
\leq\Bigl|\int_{X}\T_{1}(h_{1})\tf\cd\T_{1}(h_{2})\atf\,d\mu\Bigr|+2\eps,
}
where $\tf=S_{1}f_{1}$.
Since $\tf\in\Hwm$, by \rfr{P-cwm}, 
$\Dlim_{h}\int_{X}\T_{1}(h)\tf\cd f'\,d\mu=0$ for any $f'\in\Hwm$, 
and since $\Dl(H)>0$,
we can construct an infinite subset $B$ of $H$ such that
$\bigl|\int_{X}\T_{1}(h_{1})\tf\cd\T_{1}(h_{2})\atf\,d\mu\bigr|<\eps$ 
for any distinct $h_{1},h_{2}\in B$.
Then for any distinct $h_{1},h_{2}\in B$ we have
$\bigl|\Climsup_{n}\blan\prod_{i=1}^{k}\T_{i}(n+h_{1})f_{i},
\prod_{i=1}^{k}\T_{i}(n+h_{2})f_{i}\bran\bigr|<3\eps$,
and so by \rfr{P-vdC},
$\Climsup_{\|\cd\|,n}\prod_{i=1}^{k}\T_{i}(n)f_{i}\leq\sqrt{3\eps}$.
Since $\eps$ is arbitrary, $\Clim_{n}\prod_{i=1}^{k}\T_{i}(n)f_{i}=0$.%
\endproof
\nremark{R-jindep}{}
We defined our $\Clim${\kern1pt}s with respect to a fixed F{\o}lner sequence in $\Z$.
However, since,
for any \GLS\ $\T$ of measure preserving transformations of $X$ and any $f\in L^{2}(X)$,
$\Clim\T(n)f$ exists with respect to any F{\o}lner sequence,
this limit is the same for all F{\o}lner sequences
(since any two such sequences can be combined to produce a new one having them as subsequences).
This implies that for $\T_{1},\ld,\T_{k}\in\gT$,
the condition that $\T_{i}^{-1}\T_{j}$ for all $i\neq j$
and $\T_{1}\times\cds\times\T_{k}$ are ergodic 
is independent of the choice of a F{\o}lner sequence in $\Z$.
It now follows from \rfr{P-joint}
that the joint ergodicity of $\T_{1},\ld,\T_{k}$
does not depend on the choice of the F{\o}lner sequence either, --
which was not apriori evident.
\endnremark
For \GLSs\ based on a single transformation,
that is, of the form $\T(n)=T^{\phi(n)}$,
we have a simple criterion of ergodicity.
Recall that for a measure preserving transformation $T$ of $X$ we defined 
$$
\Eig T=\bigl\{\lam\in\C^{*}:\hbox{$Tf = \lam f$ for some $f\in L^{2}(X)$}\bigr\},
$$
and for several measure preserving transformations $T_{1},\ld,T_{k}$ of $X$,
$\Eig(T_{1},\ld,T_{k})=(\Eig T_{1})\cds(\Eig T_{k})$.
\lemma{P-sinerg}{}{}
Let $T$ be an invertible measure preserving transformation of $X$
and let $\phi$ be an unbounded \GLFs\ $\Z\ra\Z$.
Then the \GLS\ $\T(n)=T^{\phi(n)}$, $n\in\Z$, is ergodic 
iff $T$ is ergodic and $\Clim_{n}\lam^{\phi(n)}=0$ for every $\lam\in\Eig T\sm\{1\}$.
\endlemma
\proof{}
The ``only if'' part is clear.
Let $L^{2}(X)=\Hc\oplus\Hwm$ be the compact/weak mixing decomposition induced by $T$.
If $T$ is ergodic and $\Clim_{n}\lam^{\phi(n)}=0$ for every $\lam\in\Eig T\sm\{1\}$,
then $\Clim_{n}\T(n)f=\int_{X}f\,d\mu$ for any $f\in\Hc$,
so $\T$ is ergodic on $\Hc$.
Considering $T$ as a unitary operator on $\Hwm$,
immerse it into a continuous unitary flow.
Using \rfr{P-LBP}, write $\phi(n)=an+\psi(n)$, $n\in\Z$, where $a\in\R$ and $\psi$ is a \BGLF;
then $a\neq 0$, so $T^{a}$ is weakly mixing on $\Hwm$.
Since for any $f\in\Hwm$ the set $\{T^{\psi(n)}f,\ n\in\Z\}$ is totally bounded,
$\T$ is weakly mixing, and so, ergodic on $\Hwm$.%
\endproof
Here are now reincarnations of Corollaries~\rfrn{P-single} and \rfrn{P-spec}:
\corollary{P-Csingle}{}{}
Let $T$ be an invertible weakly mixing measure preserving transformation of $X$
and let $\phi_{1},\ld,\phi_{k}$ be unbounded \GLFs\ 
such that $\phi_{j}-\phi_{i}$ are unbounded for all $i\neq j$;
then the \GLSs\ $T^{\phi_{1}(n)},\ld,T^{\phi_{k}(n)}$, $n\in\Z$, are jointly ergodic.
In particular, for any distinct $\alf_{1},\ld,\alf_{k}\in\R\sm\{0\}$,
the \GLSs\ $T^{[\alf_{1}n]},\ld,T^{[\alf_{k}n]}$ are jointly ergodic.
\endcorollary
\proof{}
By \rfr{P-sinerg}, the \GLSs\ $T^{-\phi_{i}(n)}T^{\phi_{j}(n)}$ are weakly mixing
and so, ergodic for all $i\neq j$.
Reasoning the same way, 
we also see that the \GLS\ $T^{\phi_{1}(n)}\times\cds\times T^{\phi_{k}(n)}$ is weakly mixing.
By \rfr{P-joint}, $T^{\phi_{1}(n)},\ld,T^{\phi_{k}(n)}$ are jointly ergodic.%
\endproof
\corollary{P-Cspec}{}{}
Let $T_{1},\ld,T_{k}$ be commuting invertible jointly ergodic 
measure preserving transformations of $X$
and let $\phi$ be an unbounded \GLF\ $\Z\ra\Z$;
then the \GLSs\ $T_{1}^{\phi(n)},\ld,T_{k}^{\phi(n)}$, $n\in\Z$, are jointly ergodic
iff $\Clim_{n}\lam^{\phi(n)}=0$ for every $\lam\in\Eig(T_{1},\ld,T_{k})\sm\{1\}$.
In particular, for any irrational $\alf\in\R$,
the \GLSs\ $T_{1}^{[\alf n]},\ld,T_{k}^{[\alf n]}$ are jointly ergodic
iff $e^{2\pi i\alf^{-1}\Q}\cap\Eig(T_{1},\ld,T_{k})=\{1\}$.
\endcorollary
\proof{}
First of all, notice that $\Eig(T_{1},\ld,T_{k})=\Eig(T_{1}\times\cds\times T_{k})$.
Since the transformations $T_{1},\ld,T_{k}$ are ergodic, 
they share the set of eigenfunctions,
so for any $i$ and $j$ we have
$\Eig(T_{i},T_{j})\sle\Eig T_{i}\cd\Eig T_{j}\sle\Eig(T_{1},\ld,T_{k})$ as well.
Applying \rfr{P-sinerg} to the transformations $T_{i}^{-1}T_{j}$ for $i\neq j$
and $T_{1}\times\cds\times T_{k}$,
we get the first assertion.

The case $\phi(n)=[\alf n]$, with an irrational $\alf$,
is now managed by the following lemma:
\lemma{P-alfbet}{}{}
For an irrational $\alf$ and a real $\bet$ 
one has $\Clim_{n}e^{2\pi i[\alf n]\bet}=0$
iff $\alf\bet\not\in\Z\alf+\Q$.
\endlemma
\proof{}
We have $[\alf n]\bet=\alf\bet n-\{\alf n\}\bet$, $n\in\Z$.
Consider the sequence $u_{n}=(\{\alf\bet n\},\{\alf n\})$ 
in the torus $\TT^{2}_{(x,y)}=\R^{2}/\Z^{2}$,
so that the sequence $([\alf n]\bet)\mod 1$ is its image in $\TT$ 
under the mapping $\sig(x,y)=\bigl(\{x\}-\bet\{y\}\bigr)\mod 1$.
If $\alf\bet$ and $\alf$ are rationally independent modulo 1,
$(u_{n})$ is uniformly distributed in $\TT^{2}$,
thus $\sig(u_{n})$ is uniformly distributed in $\TT$,
and so, $\Clim_{n}e^{2\pi i[\alf n]\bet}=\Clim_{n}e^{2\pi i\sig(u_{n})}=0$.
Let $\alf\bet$ and $\alf$ be rationally dependent modulo 1,
$k\alf\bet=m\alf+l$, where $k\in\N$ and $m,l\in\Z$ with $\gcd(k,m,l)=1$.
Then the sequence $(u_{n})$ in uniformly distributed 
in the subgroup $S$ of $\TT^{2}$ defined by the equation $kx=my$.
$\sig$ maps $S$ to $k$ isomorphic intervals 
$\bigl[\frac{m}{k}j,\frac{m}{k}(j+1)-\bet\bigr)$, $j=0,\ld,k-1$, in $\TT$,
and the sequence $\sig(u_{n})$
is uniformly distributed in the weighted union of these intervals.
It follows that $\Clim_{n}e^{2\pi i\sig(u_{n})}=0$
unless all the intervals coincide, 
which happens iff $k\big|m$, that is, iff $\alf\bet\in\Z\alf+\Q$.
In this situation of a single interval
it is still possible that $\Clim_{n}e^{2\pi i\sig(u_{n})}=0$, --
if this interval covers $\TT$ an integer number of times,
that is, iff $\frac{m}{k}-\bet\in\Z\sm\{0\}$;
however, this is never the case since $\alf$ is irrational.
Thus, if $\alf\bet\in\Z\alf+\Q$, $\Clim_{n}e^{2\pi i\sig(u_{n})}\neq0$.%
\endproof
\frgdsp
\endproof
\section{S-Primes}{Joint ergodicity of \GLSs\ along primes}

In this section we will adapt some results from \brfr{GT} and technique from \brfr{Sun}
to establish a condition for several \GLSs\ to be jointly ergodic along primes.

By $\P$ we will denote the set of prime integers.
Let us also use the following notation:
for $N\in\N$ let $\P(N)=\P\cap\{1,\ld,N\}$, $\pi(N)=|\P(N)|$,
and $R(N)=\{r\in\{1,\ld,N\}:\gcd(r,N)=1\}$.
As above, we fix a commutative group $G$ of measure preserving transformations
of a probability measure space $(X,\mu)$
and denote by $\gT$ the group of \GLSs\ in $G$.
We will prove the following theorem:
\theorem{P-Primes}{}{}
Let $\T_{1},\ld,\T_{k}\in\gT$ be such that for any $W\in\N$ and $r\in R(W)$
the \GLSs\ $\T_{i,W,r}(n)=\T_{i}(Wn+r)$, $i=1,\ld,k$, are jointly ergodic.
Then for any $f_{1},\ld,f_{k}\in L^{\infty}(X)$,
$$
\lim_{N\ras\infty}\frac{1}{\pi(N)}\sum_{p\in\P(N)}
\T_{1}(p)f_{1}\cds\T_{k}(p)f_{k}=\prod_{i=1}^{k}\int_{X}f_{i}\,d\mu
\quad\hbox{(in $L^{2}$ norm)}.
$$
\endtheorem
\remark{}
In general, joint ergodicity of $\T_{i}$, $i=1,\ld,k$,
does not imply that of $\T_{i,W,r}$.
Indeed, let $\alf\in\R\sm\Q$, let $X=\{0,1\}$ with measure $\mu(\{0\})=\mu(\{1\})=1/2$, 
let $Tx=(x+1)\mod 2$, let $\T_{1}(n)=T^{n}$ and $\T_{2}(n)=T^{[\alf n]}$, $n\in\N$.
Then $\T_{1}$ and $\T_{2}$ are jointly ergodic, but $\T_{1}(2n+1)$ is not ergodic.
Notice also that the assertion of the theorem does not hold for these $\T_{1}$ and $\T_{2}$:
for functions $f_{1}$ and $f_{2}$ on $X$,
$\lim_{N\ras\infty}\frac{1}{\pi(N)}\sum_{p\in\P(N)}
\T_{1}(p)f_{1}\cd\T_{2}(p)f_{2}=Tf_{1}\int_{X}f_{2}\,d\mu$.
\endremark

Following \brfr{GT}, we introduce ``the modified von Mangoldt function''
$\Lam'(n)=1_{\P}(n)\log n$, $n\in\N$.
The following simple lemma allows one to rewrite the average in \rfr{P-Primes}
in terms of $\Lam'$:
\lemma{P-Lambda}{}{(Cf.\nasp\ Lemma 1 in \brfr{FHK}.)}
For any bounded sequence $(v_{n})$ of vectors in a normed vector space,
$\lim_{N\ras\infty}\bigl\|\frac{1}{\pi(N)}\sum_{p\in\P(N)}v_{p}
-\frac{1}{N}\sum_{n=1}^{N}\Lam'(n)v_{n}\bigr\|=0.$
\endlemma

A (compact) {\it nilmanifold\/} $\cN$ 
is a compact homogeneous space of a nilpotent Lie group $\cG$;
a {\it nilrotation\/} of $\cN$ is a translation by an element of $\cG$.
Nilmanifolds are characterized by the nilpotency class and the number of generators of $\cG$;
for any $k,d\in\N$ there exists a universal, ``free'' nilmanifold $\cN_{k,d}$
of nilpotency class $k$, with $d$ ``continuous'' and $d$ ``discrete'' generators%
\fnote{``A continuous generator'' of a nilmanifold $\cN$
is a continuous flow $(a^{t})_{t\in\R}$ in the group $\cG$;
``a discrete generator'' is just an element of $\cG$.}
such that any nilmanifold of class $\leq k$ and with $\leq d$ generators
is a factor of $\cN_{k,d}$.
{\it A basic nilsequence\/} is a sequence of the form $\eta(n)=g(a^{n})$
where $g$ is a continuous function on a nilmanifold $\cN$
and $a$ is a nil-rotation of $\cN$.
We may always assume that $\cN=\cN_{k,d}$ for some $k$ and $d$;
the minimal such $k$ is said to be the nilpotency class of $\eta$.
Given $k,d\in\N$ and $M>0$, 
we will denote by $\cL_{k,d,M}$ the set of basic nilsequences $\eta(n)=g(a^{n})$
where the function $g\in C(\cN_{k,d})$ is Lipschitz with constant $M$ and $|g|\leq M$.
(A smooth metric on each nilmanifold $\cN_{k,d}$ is assumed to be chosen.)

Following \brfr{GT},
for $W,r\in\N$ we define $\Lam'_{W,r}(n)=\frac{\ophi(W)}{W}\Lam'(Wn+r)$, $n\in\N$,
where $\ophi$ is the Euler function, $\ophi(W)=|R(W)|$.
By $\cW$ we will denote the set of integers of the form $W=\prod_{p\in\P(m)}p$, $m\in\N$.
It is proved in \brfr{GT} 
that ``the $W$-tricked von Mangoldt sequences $\Lam'_{W,r}$ are orthogonal to nilsequences'';
here is a weakened version of Proposition~10.2 from \brfr{GT}:
\proposition{P-vMortN}{}{}
For any $k\in\N$ and $M>0$,
$$
\lim_{\sdup{W\in\cW}{W\ras\infty}}\limsup_{N\ras\infty}
\sup_{\sdup{\eta\in\cL_{k,d,M}}{r\in R(W)}}
\Bigl|\frac{1}{N}\sum_{n=1}^{N}(\Lam'_{W,r}(n)-1)\eta(n)\Bigr|=0.
$$
\endproposition

We need to extend \rfr{P-vMortN} to sequences slightly more general than nilsequences:
\lemma{P-vMortP}{}{(Cf. \brfr{Sun}, Proposition 3.2.)}
Let $P$ be a polygonal subset of a torus $\cM$ and let $u\in\cM$.
For any $k\in\N$ and $M>0$,
$$
\lim_{\sdup{W\in\cW}{W\ras\infty}}\limsup_{N\ras\infty}
\sup_{\sdup{\eta\in\cL_{k,d,M}}{r\in R(W)}}
\Bigl|\frac{1}{N}\sum_{n=1}^{N}(\Lam'_{W,r}(n)-1)1_{P}((Wn+r)u)\eta(n)\Bigr|=0.
$$
\endlemma
\proof{}
Let $Z=\ovr{\Z u}$;
$Z$ is a finite union of subtori $\cM_{1},\ld,\cM_{l}$ of $\cM$.
Let $\eps>0$.
Choose smooth functions $g_{1},g_{2}$ on $\cM$
such that $0\leq g_{1}\leq 1_{P}\leq g_{2}\leq 1$
and the set $S=\{g_{1}\neq g_{2}\}$ is polygonal
with $\lam_{\cM_{i}}(S\cap\cM_{i})\leq\eps\lam_{\cM_{i}}(\cM_{i})$, $i=1,\ld,l$,
where $\lam_{\cM_{i}}$ is the normalized Haar measure on $\cM_{i}$.
Then for any $W$ and $r$, 
the sequences $\zeta_{1,W,r}(n)=g_{1}((Wn+r)u)$ and $\zeta_{2,W,r}(n)=g_{2}((Wn+r)u)$
are (1-step) basic nilsequences,
and since the sequence $(Wn+r)u$ is uniformly distributed 
in the union of several components of $Z$,
the set $\{n:\zeta_{1}(n)\neq\zeta_{2}(n)\}$ has density $<\eps$.
Notice that if $M'$ is the sum of $M$ and the Lipschitz's constants of $g_{1}$ and $g_{2}$
and $d'=d+\dim\cM$,
then for any $\eta\in\cL_{k,d,M}$ 
one has $\zeta_{1,W,r}\eta,\,\zeta_{2,W,r}\eta\in\cL_{k,d',M'}$.

Let $\cL_{k,d,M}^{+}=\{\eta\in\cL_{k,d,M}:\eta\geq 0\}$.
For any $W$, $r$, and any $\eta\in\cL_{k,d,M}^{+}$,
for every $n\in\N$ we have
\lequ{
(\Lam'_{W,r}(n)-1)1_{P}((Wn+r)u)\eta(n)
\leq\bigl(\Lam'_{W,r}(n)\zeta_{2,W,r}(n)-\zeta_{1,W,r}(n)\bigr)\eta(n)
\-\\\-
=(\Lam'_{W,r}(n)-1)\zeta_{2,W,r}(n)\eta(n)+\bigl(\zeta_{2,W,r}(n)-\zeta_{1,W,r}(n)\bigr)\eta(n).
}
By \rfr{P-vMortN}, $\lim_{\sdup{W\in\cW}{W\ras\infty}}\limsup_{N\ras\infty}
\sup_{\sdup{\eta\in\cL_{k,d,M}}{r\in R(W)}}
\bigl|\frac{1}{N}\sum_{n=1}^{N}(\Lam'_{W,r}(n)-1)\zeta_{2,W,r}(n)\eta(n)\bigr|=0$,
whereas for any $W$, $r$, and any $\eta\in\cL_{k,d,M}^{+}$,
$\lim_{N\ras\infty}\frac{1}{N}\sum_{n=1}^{N}
\bigl|\zeta_{2,W,r}(n)-\zeta_{1,W,r}(n)\bigr|\eta(n)\leq M\eps$;
thus,
$$
\limsup_{\sdup{W\in\cW}{W\ras\infty}}\limsup_{N\ras\infty}
\sup_{\sdup{\eta\in\cL_{k,d,M}^{+}}{r\in R(W)}}
\frac{1}{N}\sum_{n=1}^{N}(\Lam'_{W,r}(n)-1)1_{P}((Wn+r)u)\eta(n)\leq M\eps.
$$
Similarly,
\equ{
(\Lam'_{W,r}(n)-1)1_{P}((Wn+r)u)\eta(n)
\geq(\Lam'_{W,r}(n)-1)\zeta_{1,W,r}(n)\eta(n)-\bigl(\zeta_{2,W,r}(n)-\zeta_{1,W,r}(n)\bigr)\eta(n),
}
so
$$
\liminf_{\sdup{W\in\cW}{W\ras\infty}}\liminf_{N\ras\infty}
\inf_{\sdup{\eta\in\cL_{k,d,M}^{+}}{r\in R(W)}}
\frac{1}{N}\sum_{n=1}^{N}(\Lam'_{W,r}(n)-1)1_{P}((Wn+r)u)\eta(n)\geq-M\eps.
$$
Hence,
$$
\lim_{\sdup{W\in\cW}{W\ras\infty}}\limsup_{N\ras\infty}
\sup_{\sdup{\eta\in\cL_{k,d,M}^{+}}{r\in R(W)}}
\Bigl|\frac{1}{N}\sum_{n=1}^{N}(\Lam'_{W,r}(n)-1)1_{P}((Wn+r)u)\eta(n)\Bigr|=0;
$$
since $\cL_{k,d,M}=\cL_{k,d,M}^{+}-\cL_{k,d,M}^{+}$, we are done.%
\endproof

Let $k,N\in\N$; 
for sequences $b\col\{1,\ld,N\}\ra\R$ we define {\it the $k$-th Gowers's norm\/} by
$$
\|b\|_{U^{k}[N]}=\Bigl(\frac{1}{N^{k}}\sum_{h_{1},\ld,h_{k}=1}^{N}
\Bigl|\frac{1}{N}\sum_{n=1}^{N-h_{1}-\cds-h_{k}}
\prod_{e_{1},\ld,e_{k}\in\{0,1\}}b(n+e_{1}h_{1}+\cds+e_{k}h_{k})\Bigr|\Bigr)^{1/2^{k}}
$$
(where we assume $\sum_{n=1}^{m}=0$ if $m\leq 0$).
The next result we need is the fact that, on a certain class of sequences,
``the $k$-th Gowers norm is continuous with respect to the system of seminorms 
$\|b\|_{\eta}=\bigl|\frac{1}{N}\sum_{n=1}^{N}b(n)\eta(n)\bigr|$, $\eta\in\cL_{k,d,M}$''.
To avoid unnecessary technicalities,
we will only formulate the following lemma,
which is a corollary of Propositions~10.1 and 6.4 in \brfr{GT}:
\lemma{P-avGow}{}{}
For any $\eps>0$ and $k\in\N$ there exist $d\in\N$, $M>0$, and $\del>0$ 
such that for any $N\in\N$,
if a sequence $b\col\{1,\ld,N\}\ra\R$ 
satisfies $|b|\leq 1+\Lam'_{W,r}$ for some $W\in\cW$ and $r\in R(W)$
and $\sup_{\eta\in\cL_{k,d,M}}\bigl|\frac{1}{N}\sum_{n=1}^{N}b(n)\eta(n)\bigr|<\del$,
then $\|b\|_{U^{k}[N]}<\eps$.
\endlemma
\remark{}
Proposition~10.1 was proved in \brfr{GT} modulo the ``Inverse Gowers-norm Conjecture'',
which has then been confirmed in \brfr{GTZ}.
\endremark
Combining \rfr{P-vMortP} and \rfr{P-avGow}, 
applied to the sequence $b(n)=(\Lam'_{W,r}(n)-1)1_{P}((Wn+r)u)$, we obtain:
\lemma{P-Gow}{}{(Cf.\nasp\ \brfr{Sun}, Proposition~3.2.)}
Let $P$ be a polygonal region in a torus $\cM$ and let $u\in\cM$.
Then for any $k\in\N$,
$$
\lim_{\sdup{W\in\cW}{W\ras\infty}}\limsup_{N\ras\infty}\max_{r\in R(W)}
\bigl\|(\Lam'_{W,r}(n)-1)1_{P}((Wn+r)u)\bigr\|_{U^{k}[N]}=0.
$$
\endlemma
From \rfr{P-Gow} we now deduce:
\proposition{P-Lamav}{}{(Cf.\nasp\ \brfr{Sun}, Proposition 4.1)}
Let $\T_{1},\ld,\T_{k}\in\gT$ and $f_{1},\ld,f_{k}\in L^{\infty}(X)$.
For any $f_{1},\ld,f_{k}\in L^{\infty}(X)$ we have
$$
\lim_{\sdup{W\in\cW}{W\ras\infty}}\limsup_{N\ras\infty}\max_{r\in R(W)}
\Bigl\|\frac{1}{N}\sum_{n=1}^{N}(\Lam'_{W,r}(n)-1)
\T_{1}(Wn+r)f_{1}\cds\T_{k}(Wn+r)f_{k}\Bigr\|_{L^{2}(X)}=0.
$$
\endcorollary
\proof{}
We will assume that $|f_{i}|\leq 1$, $i=1,\ld,k$.
Let, by \rfr{P-perT}(iii), 
$\cM$ be a torus, $u\in\cM$,
and $\cM=\bigcup_{j=1}^{l}P_{j}$ be a polygonal partition of $\cM$
such that for every $i$ and $j$,
$\T_{i}(n)^{-1}\T_{i}(n+h)$ does not depend on $n$ whenever both $nu,(n+h)u\in P_{j}$.
We will show that for any $j$ and any $W,r\in\N$,
\lequ{
\Bigl\|\frac{1}{N}\sum_{n=1}^{N}(\Lam'_{W,r}(n)-1)1_{P_{j}}((Wn+r)u)
\T_{1}(Wn+r)f_{1}\cds\T_{k}(Wn+r)f_{k}\Bigr\|_{L^{2}(X)}
\-\\\-
\leq 2\bigl\|(\Lam'_{W,r}(n)-1)1_{P_{j}}((Wn+r)u)\bigr\|_{U^{k}[N]}+o_{N}(1);
}
via \rfr{P-Gow}, this will imply that
$$
\lim_{\sdup{W\in\cW}{W\ras\infty}}\limsup_{N\ras\infty}\max_{r\in R(W)}
\Bigl\|\frac{1}{N}\sum_{n=1}^{N}(\Lam'_{W,r}(n)-1)1_{P_{j}}((Wn+r)u)
\T_{1}(Wn+r)f_{1}\cds\T_{k}(Wn+r)f_{k}\Bigr\|_{L^{2}(X)}=0
$$
for each $j=1,\ld,l$, from which the assertion of the proposition follows.

Fix $j$, $W$, and $r$;
put $P=P_{j}$, 
$b(n)=(\Lam'_{W,r}(n)-1)1_{P}((Wn+r)u)$ for $n\in\N$,
and $\tT_{i}(n)=\T_{i}(Wn+r)$ for $i=1,\ld,k$.

By \rfr{P-vdC2}, for any $N$,
\lequ{
\Bigl\|\frac{1}{N}\sum_{n=1}^{N}b(n)
\tT_{1}(n)f_{1}\cds\tT_{k}(n)f_{k}\Bigr\|_{L^{2}(X)}^{2}
\-\\\-
\leq\frac{2}{N}\sum_{h_{1}=1}^{N}\Bigl|\frac{1}{N}\sum_{n=1}^{N-h_{1}}
\int_{X}b(n)b(n+h_{1})\cd\tT_{1}(n)f_{1}\cd\tT_{1}(n+h_{1})\af_{1}\cds
\tT_{k}(n)f_{k}\cd\tT_{k}(n+h_{1})\af_{k}\,d\mu\Bigr|
\-\\\-
+\frac{1}{N^{2}}\sum_{n=1}^{N}|b(n)|^{2}
\|f_{1}\|_{L^{\infty}(X)}^{2}\cds\|f_{k}\|_{L^{\infty}(X)}^{2}.
}
Since $|b(n)|\leq\log(Wn+r)$ and $|f_{1}|,\ld,|f_{k}|\leq 1$,
the second summand is $o(1)$ as $N\ras\infty$.
By the definition of $P=P_{j}$, 
for each $i$ there exists a sequence $\S_{i}(h_{1})$, $h_{1}\in\Z$, of transformations
such that
$\tT_{i}(n)^{-1}\tT_{i}(n+h_{1})=\S_{i}(h_{1})$ if $1_{P}((Wn+r)u)1_{P}((W(n+h_{1})+r)u)\neq 0$,
and so, if $b(n)b(n+h_{1})\neq 0$.
Thus, if we put 
$f_{i,h_{1}}=f_{i}\cd\S_{i}(h_{1})\af_{i}$, $i=1,\ld,k$, $h_{1}\in\N$,
we get
\lequ{
\Bigl\|\frac{1}{N}\sum_{n=1}^{N}b(n)
\tT_{1}(n)f_{1}\cds\tT_{k}(n)f_{k}\Bigr\|_{L^{2}(X)}^{2}
\-\\\kern30mm
\leq\frac{2}{N}\sum_{h_{1}=1}^{N}\Bigl|\frac{1}{N}\sum_{n=1}^{N-h_{1}}
\int_{X}b(n)b(n+h_{1})\cd\tT_{1}(n)f_{1,h_{1}}\cds\tT_{k}(n)f_{k,h_{1}}\,d\mu\Bigr|
+o(1)
\-\\\kern30mm
=\frac{2}{N}\sum_{h_{1}=1}^{N}
\Bigl|\int_{X}f_{1,h_{1}}\frac{1}{N}\sum_{n=1}^{N-h_{1}}b(n)b(n+h_{1})\cd
(\tT_{1}^{-1}\tT_{2})(n)f_{2,h_{1}}\cds(\tT_{1}^{-1}\tT_{k})(n)f_{k,h_{1}}\,d\mu\Bigr|
+o(1)
\-\\\kern30mm
\leq\frac{2}{N}\sum_{h_{1}=1}^{N}
\Bigl\|\frac{1}{N}\sum_{n=1}^{N-h_{1}}b(n)b(n+h_{1})\cd
(\tT_{1}^{-1}\tT_{2})(n)f_{2,h_{1}}\cds(\tT_{1}^{-1}\tT_{k})(n)f_{k,h_{1}}\Bigr\|_{L^{2}(X)}
+o(1).
\-}
In the same way, for every $h_{1}$,
\lequ{
\Bigl\|\frac{1}{N}\sum_{n=1}^{N-h_{1}}b(n)b(n+h_{1})\cd
(\tT_{1}^{-1}\tT_{2})(n)f_{2,h_{1}}\cds(\tT_{1}^{-1}\tT_{k})(n)f_{k,h_{1}}
\Bigr\|_{L^{2}(X)}^{2}
\-\\\-
\leq\frac{2}{N}\sum_{h_{2}=1}^{N-h_{1}}
\Bigl\|\frac{1}{N}\sum_{n=1}^{N-h_{1}-h_{2}}b(n)b(n+h_{1})b(n+h_{2})b(n+h_{1}+h_{2})
\-\\\-
\cd(\tT_{2}^{-1}\tT_{3})(n)f_{3,h_{1},h_{2}}\cds(\tT_{2}^{-1}\tT_{k})(n)f_{k,h_{1},h_{2}}
\Bigr\|_{L^{2}(X)}+o(1),
}
for some functions $f_{i,h_{1},h_{2}}$ of modulus $\leq1$, 
and so, by Schwarz's inequality,
\lequ{
\Bigl\|\frac{1}{N}\sum_{n=1}^{N}b(n)
\tT_{1}(n)f_{1}\cds\tT_{k}(n)f_{k}\Bigr\|_{L^{2}(X)}^{4}
\-\\\-
\leq\frac{2^{3}}{N^{2}}\sum_{h_{1},h_{2}=1}^{N}
\Bigl\|\frac{1}{N}\sum_{n=1}^{N-h_{1}-h_{2}}b(n)b(n+h_{1})b(n+h_{2})b(n+h_{1}+h_{2})
\-\\\-
\cd(\tT_{2}^{-1}\tT_{3})(n)f_{3,h_{1},h_{2}}\cds(\tT_{2}^{-1}\tT_{k})(n)f_{k,h_{1},h_{2}}
\Bigr\|_{L^{2}(X)}
+o(1).
}
(We always assume that $\sum_{n=1}^{m}=0$ if $m\leq 0$.)
Applying \rfr{P-vdC2} $k-2$ more times, we arrive at
\lequ{
\Bigl\|\frac{1}{N}\sum_{n=1}^{N}b(n)
\tT_{1}(n)f_{1}\cds\tT_{k}(n)f_{k}\Bigr\|_{L^{2}(X)}^{2^{k}}
\-\\\-
\leq\frac{2^{2^{k}-1}}{N^{k}}\sum_{h_{1},\ld,h_{k}=1}^{N}
\Bigl|\frac{1}{N}\sum_{n=1}^{N-h_{1}-\cds-h_{k}}
\prod_{e_{1},\ld,e_{k}\in\{0,1\}}b(n+e_{1}h_{1}+\cds+e_{k}h_{k})\Bigr|+o(1)
\-\\\-
=2^{2^{k}-1}\|b\|_{U^{k}[N]}^{2^{k}}+o(1).
}
\frgdsp\endproof
We are now in position to prove \rfr{P-Primes}:
\proof{of \rfr{P-Primes}}
For short, put $\ff(n)=\T_{1}(n)f_{1}\cds\T_{k}(n)f_{k}$, $n\in\N$.
By \rfr{P-Lambda}, we have to show that
$\lim_{N\ras\infty}\frac{1}{N}\sum_{n=1}^{N}\Lam'(n)\ff(n)=\prod_{i=1}^{k}\int_{X}f_{i}\,d\mu$.
Let $\eps>0$.
By \rfr{P-Lamav}, we can choose $W\in\cW$ such that
for any $N$ large enough and any $r\in R(W)$ one has
$$
\Bigl\|\frac{1}{N}\sum_{n=1}^{N}(\Lam'_{W,r}(n)-1)\ff(Wn+r)\Bigr\|_{L^{2}(X)}<\eps,
$$
and so,
$$
\Bigl\|\frac{1}{NW}\sum_{n=1}^{N}\Lam'(Wn+r)\ff(Wn+r)
-\frac{1}{N\ophi(W)}\sum_{n=1}^{N}\ff(Wn+r)\Bigr\|_{L^{2}(X)}<\frac{\eps}{\ophi(W)}.
$$
Summing this up for all $r\in R(W)$, 
and taking into account that $\Lam'(Wn+r)=0$ if $r\not\in R(W)$,
we obtain
$$
\Bigl\|\frac{1}{NW}\sum_{n=1}^{NW}\Lam'(n)\ff(n)
-\frac{1}{\ophi(W)}\sum_{r\in R(W)}\frac{1}{N}\sum_{n=1}^{N}\ff(Wn+r)\Bigr\|_{L^{2}(X)}<\eps.
$$
By the theorem's assumption, for any $r\in R(W)$,
$\bigl\|\frac{1}{N}\sum_{n=1}^{N}\ff(Wn+r)
-\prod_{i=1}^{k}\int_{X}f_{i}\,d\mu\bigr\|_{L^{2}(X)}<\eps$
for all $N$ large enough.
Hence, $\bigl\|\frac{1}{NW}\sum_{n=1}^{NW}\Lam'(n)\ff(n)\bigr\|_{L^{2}(X)}<2\eps$ for such $N$,
and so, 
$$
\lim_{N\ras\infty}\frac{1}{N}\sum_{n=1}^{N}\Lam'(n)\ff(n)
=\lim_{N\ras\infty}\frac{1}{NW}\sum_{n=1}^{NW}\Lam'(n)\ff(n)
=\prod_{i=1}^{k}\int_{X}f_{i}\,d\mu.
$$
\frgdsp\endproof

We will now collect some special cases of \rfr{P-Primes}.
It was shown in \brfr{Berend} that if $T_{1},\ld,T_{k}$, with $k\geq 2$,
are commuting, invertible, jointly ergodic measure preserving transformations,
then they are actually totally jointly ergodic,
that is, for any $W\in\N$ and $r\in\Z$,
$T_{1}^{Wn+r},\ld,T_{k}^{Wn+r}$ are jointly ergodic.
Hence, by \rfr{P-Primes}, we obtain:
\theorem{P-totjerg}{}{}
Let $T_{1},\ld,T_{k}$, where $k\geq 2$,
be commuting, invertible, jointly ergodic measure preserving transformations of $X$.
Then for any $f_{1},\ld,f_{k}\in L^{\infty}(X)$, in the $L^{2}$-norm,
$$
\lim_{N\ras\infty}\frac{1}{\pi(N)}\sum_{p\in\P(N)}
T_{1}^{p}f_{1}\cds T_{k}^{p}f_{k}=\prod_{i=1}^{k}\int_{X}f_{i}\,d\mu.
$$
\endtheorem
The following is a corollary of \rfr{P-Primes} and \rfr{P-single}:
\corollary{P-psingle}{}{}
Let $T$ be a weakly mixing invertible measure preserving transformation of $X$
and let $\phi_{1},\ld,\phi_{k}$ be unbounded \GLFs\ $\Z\ra\Z$
such that $\phi_{j}-\phi_{i}$ are unbounded for all $i\neq j$. 
Then for any $f_{1},\ld,f_{k}\in L^{\infty}(X)$,
$$
\lim_{N\ras\infty}\frac{1}{\pi(N)}\sum_{p\in\P(N)}T^{\phi_{1}(p)}f_{1}\cds T^{\phi_{k}(p)}f_{k}
=\prod_{i=1}^{k}\int_{X}f_{i}\,d\mu.
$$
In particular, for any distinct $\alf_{1},\ld,\alf_{k}\in\R\sm\{0\}$,
$$
\lim_{N\ras\infty}\frac{1}{\pi(N)}\sum_{p\in\P(N)}T^{[\alf_{1}p]}f_{1}\cds T^{[\alf_{k}p]}f_{k}
=\prod_{i=1}^{k}\int f_{i}\,d\mu.
$$
\endcorollary
From \rfr{P-Primes} and \rfr{P-spec} we obtain:
\corollary{P-pspec}{}{}
Let $T_{1},\ld,T_{k}$ be commuting invertible jointly ergodic 
measure preserving transformations of $X$
and let $\phi$ be an unbounded \GLF\ $\Z\ra\Z$
such that $\lim_{N\ras\infty}\frac{1}{N}\sum_{n=1}^{N}\lam^{\phi(Wn+r)}=0$
for every $\lam\in\Eig(T_{1},\ld,T_{k})$, $W\in\cW$, and $r\in R(W)$.
Then for any $f_{1},\ld,f_{k}\in L^{\infty}(X)$,
$$
\lim_{N\ras\infty}\frac{1}{\pi(N)}\sum_{p\in\P(N)}T_{1}^{\phi(p)}f_{1}\cds T_{k}^{\phi(p)}f_{k}
=\prod_{i=1}^{k}\int f_{i}\,d\mu\
$$
In particular, if $\alf\in\R$ is irrational 
and such that $e^{2\pi i\alf^{-1}\Q}\cap\Eig(T_{1},\ld,T_{k})=\{1\}$,
then
$$
\lim_{N\ras\infty}\frac{1}{\pi(N)}\sum_{p\in\P(N)}T_{1}^{[\alf p]}f_{1}\cds T_{k}^{[\alf p]}f_{k}
=\prod_{i=1}^{k}\int f_{i}\,d\mu\
\hbox{for any $f_{1},\ld,f_{k}\in L^{\infty}(X)$}.
$$ 
\endcorollary
\section{S-cont}{\GLFas\ of a continuous parameter}

Let $\T(t)$, $t\in\R$, be a family of measure preserving transformations of $X$.
We say that $\T$ is {\it ergodic\/} if, for any $f\in L^{2}(X)$,
$\lim_{b\ras\infty}\frac{1}{b}\int_{0}^{b}\T(t)f\,dt=\int_{X}f\,d\mu$ in $L^{2}$-norm,
and {\it uniformly ergodic\/} if, for any $f\in L^{2}(X)$,
$\lim_{b\ras\infty}\frac{1}{b-a}\int_{a}^{b}\T(t)f\,dt=\int_{X}f\,d\mu$.
Given several families $\T_{1}(t),\ld,\T_{k}(t)$, $t\in\R$,
of measure preserving transformations of $X$,
we say that $\T_{1},\ld,\T_{k}$ are {\it jointly ergodic\/} if 
$$
\lim_{b\ras\infty}\frac{1}{b}\int_{0}^{b}\T_{1}(t)f_{1}\cds\T_{k}(t)f_{k}\,dt
=\prod_{i=1}^{k}\int_{X}f_{i}\,d\mu
$$ 
in $L^{2}$-norm for any $f_{1},\ld,f_{k}\in L^{\infty}(X)$,
and {\it uniformly jointly ergodic\/} if 
$$
\lim_{b-a\ras\infty}\frac{1}{b-a}\int_{a}^{b}\T_{1}(t)f_{1}\cds\T_{k}(t)f_{k}\,dt
=\prod_{i=1}^{k}\int_{X}f_{i}\,d\mu
$$ 
for any $f_{1},\ld,f_{k}\in L^{\infty}(X)$.

Let $G$ be a commutative group of measure preserving transformations of $X$.
In analogy with the terminology adopted in previous sections,
we will call a family $\T(t)$, $t\in\R$, of elements of $G$
{\it a \GLFa\/} if it is of the form
$\T(t)=T_{1}^{\phi_{1}(t)}\cds T_{r}^{\phi_{r}(t)}$, $t\in\R$,
where $T_{1},\ld,T_{r}$ are continuous homomorphisms $\R\ra G$
and $\phi_{1},\ld,\phi_{r}$ are \GLFs\ $\R\ra\R$.
Let $\gTR$ denote the set of \GLFas\ of transformations from $G$.
We have the following analogue of \rfr{P-joint}:
\theorem{P-Rjoint}{}{}
Let $\T_{1},\ld,\T_{k}\in\gTR$.
Then the following are equivalent:\\
{\rm(i)} $\T_{1},\ld,\T_{k}$ are jointly ergodic;\\
{\rm(ii)} $\T_{1},\ld,\T_{k}$ are uniformly jointly ergodic;\\
{\rm(iii)} the \GLFas\ $\T_{i}^{-1}\T_{j}$ are ergodic for all $i\neq j$
and the \GLFa\ $\T_{1}\times\cds\times\T_{k}$ is ergodic.
\endtheorem

One can verify that a (properly modified) proof of \rfr{P-joint}
works in the situation at hand as well.
An alternative and simpler approach is to derive \rfr{P-Rjoint} from \rfr{P-joint}
with the help of the techniques developed in \brfr{dcc}.
Namely, we can use the following fact:
\theorem{P-dcc}{}{(\brfr{dcc})}
Let $\tau\col\R\ra V$ be a bounded measurable mapping to a Banach space $V$
such that for every $t\in\R$, 
the limit $L_{t}=\lim_{N\ras\infty}\frac{1}{N}\sum_{n=1}^{N}\tau(nt)$
(respectively, $L_{t}=\lim_{N-M\ras\infty}\frac{1}{N-M}\sum_{n=M+1}^{N}\tau(nt)$)
exists for \ae $t\in\R$.
Then the limit $L=\lim_{b\ras\infty}\frac{1}{b}\int_{0}^{b}\tau(t)\,dt$ 
(respectively, $L=\lim_{b\ras\infty}\frac{1}{b-a}\int_{a}^{b}\tau(t)\,dt$)
also exists,
and $L_{t}=L$ for \ae $t\in\R$.
\endtheorem

To apply this result we need to verify that for any \GLFa\ $\T$ and any $t\in\R$
the sequence $\T(nt)$, $n\in\Z$, is a \GLS.
This is indeed so:
any \GLF\ $\phi$ can be written in the form $\phi(t)=\sum_{j=1}^{l}[\phi_{j}(t)]a_{j}+ct+d$,
where $\phi_{j}$ are \GLFs\ and $a_{j},c,d\in\R$,
thus for any flow $T$ and any $t\in\R$,
$T^{\phi(nt)}=\bigl(\prod_{j=1}^{l}(T^{a_{j}})^{[\phi_{j}(nt)]}\bigr)(T^{ct})^{n}T^{d}$,
in which expression all the factors are \GLSs\ in the group
generated by the transformations $T^{a_{1}},\ld,T^{a_{l}},T^{ct},T^{d}$.
We may now apply \rfr{P-dcc} in conjunction with \rfr{P-clim} 
to the family $\tau(t)=\T(t)f$,
where $\T$ is a \GLFa\ and $f\in L^{2}(X)$,
and see that the limits $\lim_{b\ras\infty}\frac{1}{b}\int_{0}^{b}\T(t)f\,dt$
and $\lim_{b-a\ras\infty}\frac{1}{b-a}\int_{a}^{b}\T(t)f\,dt$ exist,
and $\T$ is ergodic and is uniformly ergodic 
iff the \GLSs\ $\T(nt)$, $n\in\Z$, are ergodic (= uniformly ergodic) for almost all $t\in\R$.

\proof{of \rfr{P-Rjoint}}
Assume that the \GLFas\ $\T_{1},\ld,\T_{k}$ are jointly ergodic.
For any distinct $i$ and $j$, $\T_{i}$ and $\T_{j}$ are jointly ergodic,
which implies that $\T_{i}^{-1}\T_{j}$ is ergodic.
It remains to show that $\T_{1}\times\cds\times\T_{k}$ is ergodic.
$\T_{1},\ld,\T_{k}$ are ergodic, 
so the \GLSs\ $\T_{1}(nt),\ld,\T_{k}(nt)$, $n\in\Z$, are ergodic for \ae $t\in\R$.
Thus, by \rfr{P-tens},
$\Clim_{n}\bigotimes_{i=1}^{k}\T_{i}(nt)f_{i}=0$ for \ae $t\in\R$ 
whenever $f_{i}\in\Hwm$ for at least one of $i$;
by \rfr{P-dcc}, this implies that
$\lim_{b\ras\infty}\frac{1}{b}\int_{0}^{b}\prod_{i=1}^{k}\T_{i}(t)f_{i}\,dt=0$
whenever $f_{i}\in\Hwm$ for at least one of $i$.
If $f_{i}\in\Hc$ for all $i$,
we may assume that all $f_{i}$ are nonconstant eigenfunctions of elements of $G$,
so $\T_{i}(t)f_{i}=\lam_{i}(t)f_{i}$, $i=1,\ld,k$, 
where $\lam_{i}$ are functions $\R\ra\C$.
In this case,
$$
0=\lim_{b\ras\infty}\frac{1}{b}\int_{0}^{b}\prod_{i=1}^{k}\T_{i}(t)f_{i}\,dt
=\Bigl(\lim_{b\ras\infty}\frac{1}{b}\int_{0}^{b}\prod_{i=1}^{k}\lam_{i}(t)\,dt\Bigr)
\prod_{i=1}^{k}f_{i},
$$
so $\lim_{b\ras\infty}\frac{1}{b}\int_{0}^{b}\prod_{i=1}^{k}\lam_{i}(t)\,dt=0$,
so
$$
\lim_{b\ras\infty}\frac{1}{b}\int_{0}^{b}\bigotimes_{i=1}^{k}\T_{i}(t)f_{i}\,dt
=\Bigl(\lim_{b\ras\infty}\frac{1}{b}\int_{0}^{b}\prod_{i=1}^{k}\lam_{i}(t)\,dt\Bigr)
\bigotimes_{i=1}^{k}f_{i}=0.
$$
Hence, $\T_{1}\times\cds\times\T_{k}$ is ergodic.

Conversely,
if $\T_{i}^{-1}\T_{j}$ for all $i\neq j$ and $\T_{1}\times\cds\times\T_{k}$ are ergodic,
then by \rfr{P-clim} and \rfr{P-dcc},
the \GLSs\ $(\T_{i}^{-1}\T_{j})(nt)$ and $(\T_{1}\times\cds\times\T_{k})(nt)$ 
are ergodic for \ae $t\in\R$,
so, by \rfr{P-joint},
the \GLSs\ $\T_{1}(nt),\ld,\T_{k}(nt)$ are jointly ergodic 
(= uniformly jointly ergodic, see remark after the proof of \rfr{P-joint}) for \ae $t\in\R$,
so $\T_{1},\ld,\T_{k}$ are jointly ergodic and uniformly jointly ergodic by \rfr{P-dcc}.%
\endproof
For a continuous parameter,
Corollaries~\rfrn{P-single} and \rfrn{P-spec} take the following form:
\corollary{P-Rsingle}{}{}
Let $T^{s}$, $s\in\R$, be a weakly mixing continuous flow 
of measure preserving transformations of $X$
and let $\phi_{1},\ld,\phi_{k}$ be unbounded \GLFs\ $\R\ra\R$
such that $\phi_{j}-\phi_{i}$ are unbounded for all $i\neq j$;
then the \GLFas\ $T^{\phi_{1}(t)},\ld,T^{\phi_{k}(t)}$, $t\in\R$, are jointly ergodic.
In particular, for any distinct $\alf_{1},\ld,\alf_{k}\in\R\sm\{0\}$,
the \GLFas\ $T^{[\alf_{1}t]},\ld,T^{[\alf_{k}t]}$, $t\in\R$, are jointly ergodic.
\endcorollary
\corollary{P-Rspec}{}{}
Let $T_{1}^{s},\ld,T_{k}^{s}$, $s\in\R$, 
be commuting jointly ergodic continuous flows of measure preserving transformations of $X$
and let $\phi$ be an unbounded \GLF;
then the \GLFas\ $T_{1}^{\phi(t)},\ld,T_{k}^{\phi(t)}$, $t\in\R$, are jointly ergodic
iff $\lim_{b\ras\infty}\frac{1}{b}\int_{0}^{b}\lam^{\phi(t)}\,dt=0$ 
for every $\lam\in\Eig(T_{1}^{1},\ld,T_{k}^{1})\sm\{1\}$.
\endcorollary

We would also like to remark that, in the case of continuous parameter,
by using a ``change of variable'' trick
one can easily extend the results above, proved for \GLFas, 
to more general families of transformations of the form $\T(\sig(t))$, 
where $\T$ is a \GLFa\ and $\sig$ is a monotone function of ``regular'' growth.
What we mean is the following proposition:
\proposition{P-Rpow}{}{}
Let $\T_{1}(t),\ld,\T_{k}(t)$, $t\in\R$, be jointly ergodic families
of measure preserving transformations of $X$
and let $\sig\col\R\ra\R$ be a strictly increasing $C^{1}$-function
such that $\sig'$ is monotone and
$$
\lim_{b-a\ras\infty}\frac{(\sig^{-1})'(a)}{\sig^{-1}(b)-\sig^{-1}(a)}
=\lim_{b-a\ras\infty}\frac{(\sig^{-1})'(b)}{\sig^{-1}(b)-\sig^{-1}(a)}=0.
$$
Then the families $\T_{1}(\sig(t)),\ld,\T_{k}(\sig(t))$ are also jointly ergodic.
\endproposition

\remark{}
Of course, \rfr{P-Rpow} remains true when the families $\T_{i}$
are only defined on a ray $[r,\infty)$.
In this form,
it applies when $\sig$ is of the form $\sig(t)=\sum_{i=1}^{d}c_{i}t^{\alf_{i}}$,
where $\alf_{i}$ are nonnegative reals,
on the ray $[r,\infty)$ where $\sig'$ becomes monotone.
Moreover, for $\sig$ of this sort,
$\sig^{-1}$ also satisfies the assumptions of the proposition,
thus $\T_{1},\ld,\T_{k}$ are jointly ergodic 
iff $\T_{1}(\sig(t)),\ld,\T_{k}(\sig(t))$ are.
\endremark

Let us say that a function $\tau\col\R\ra\R$ {\it has uniform Ces\`aro limit $L$}
if $\lim_{b-a\ras\infty}\dsc\frac{1}{b-a}\int_{a}^{b}\tau(t)\,dt=L$.
\rfr{P-Rpow} is simply a special case of the following general fact:
\proposition{P-subs}{}{}
Let $\sig$ be as in \rfr{P-Rpow}.
Then, if a bounded function $\tau\col\R\ra\R$ has uniform Ces\`aro limit $L$,
the function $\tau(\sig(t))$ also does.
\endproposition
This proposition must be well known to aficionados,
but since we have not been able to find any references,
we will sketch its proof.
Making the substitution $s=\sig(t)$ we get
$$
\lim_{b-a\ras\infty}\frac{1}{b-a}\int_{a}^{b}\tau(\sig(t))\,dt
=\lim_{q-p\ras\infty}\frac{1}{\sig^{-1}(q)-\sig^{-1}(p)}
\int_{p}^{q}\tau(s)(\sig^{-1})'(s)ds.
$$
What we have in the right hand part of this formula,
$\lim_{b-a\ras\infty}\frac{1}{\int_{a}^{b}\om}\int_{a}^{b}\tau(t)\om(t)\,dt$,
is {\it the weighted uniform Ces\`aro limit of $\tau$ with weight $\om=(\sig^{-1})'$}.
Rewriting \rfr{P-subs} in terms of $\om$,
we reduce it to the following lemma:
\lemma{P-waver}{}{}
Let $\om\col\R\ra\R$ be a positive monotone function 
with the property that, for any $c>0$,
$\lim_{b-a\ras\infty}\om(a)/\int_{a}^{b}\om\dsc
=\lim_{b-a\ras\infty}\om(b)/\int_{a}^{b}\om=0$.
Then, if a bounded function $\tau\col[0,\infty)\ra\R$ has uniform Ces\`aro limit $L$,
then the weighted uniform Ces\`aro limit of $\tau$ with weight $\om$ is equal to $L$
(and, in particular, exists).
\endlemma
\proof{}
We will assume that $\om$ is increasing,
the case of decreasing $\om$ is similar.
\ignore
First of all, notice that for any $c>0$, the quotient $\om(x+c)/\om(x)$ is bounded.
Indeed, assume it is not.
Find $z>0$ such that $\om(b)/\int_{a}^{b}\om<\frac{1}{2c}$ whenever $b-a>z$.
Find $b$ such that $\om(b+c)=r\om(b)$, where $r=\frac{b-a}{c}$,
and put $a=b-z$.
Then we have
$$
\frac{\om(b+c)}{\int_{a}^{b+c}\om}
\geq\frac{\om(b+c)}{\om(b)(b-a)+\om(b+c)c}
=\frac{r\om(b)}{\om(b)(b-a)+r\om(b)c}
=\frac{r}{b-a+rc}=\frac{1}{2c},
$$
which contradicts the choice of $z$.
It follows that for any $c>0$, $\lim_{b-a\ras\infty}\om(b+c)/\int_{a}^{b}=0$ as well.
\endignore
Let $M=\sup|\tau|$.
Let $\eps>0$.
Find $c>0$ such that $\frac{1}{c}\int_{x}^{x+c}\tau(t)\,dt\app^{\eps}L$ for every $x>0$.
Averaging this equation with weight $\om$ over an interval $[a,b]$ 
and changing the order of integration, we get
\equ{
L\app^{\eps}
\frac{1}{\int_{a}^{b}\om}
\int_{a}^{b}\om(x)\Bigl(\frac{1}{c}\int_{x}^{x+c}\tau(t)\,dt\Bigr)\,dx
=\frac{1}{c\int_{a}^{b}\om}\int_{a}^{b}\Bigl(\int_{t-c}^{t}\om(x)\,dx\Bigr)\tau(t)\,dt
\-\\\-\kern30mm
-\frac{1}{c\int_{a}^{b}\om}\int_{a}^{a+c}\Bigl(\int_{t-c}^{a}\om(x)\,dx\Bigr)\tau(t)\,dt
+\frac{1}{c\int_{a}^{b}\om}\int_{b}^{b+c}\Bigl(\int_{t-c}^{b}\om(x)\,dx\Bigr)\tau(t)\,dt.
}
The moduli of the second and of the third summands in the right hand part of this epuality
are majorized by $\frac{M\om(b)c^{2}/2}{c\int_{a}^{b}\om}$
and tend to 0 as $b-a\ra\infty$.
We now claim that, for large $b-a$, 
the first summand is close to
$\frac{1}{\int_{a}^{b}\om}\int_{a}^{b}\tau(t)\om(t)\,dt$.
Indeed, taking into account the monotonicity of $\om$, we have
\lequ{
\left|\frac{1}{c\int_{a}^{b}\om}
\int_{a}^{b}\Bigl(\int_{t-c}^{t}\om(x)\,dx\Bigr)\tau(t)\,dt
-\frac{1}{\int_{a}^{b}\om}\int_{a}^{b}\tau(t)\om(t)\,dt\right|
=\frac{1}{c\int_{a}^{b}\om}
\left|\int_{a}^{b}\Bigl(\int_{t-c}^{t}\bigl(\om(x)-\om(t)\bigr)\,dx\Bigr)\tau(t)\,dt\right|
\-\\\-
\leq\frac{M}{c\int_{a}^{b}\om}
\int_{a}^{b}\Bigl(\int_{t-c}^{t}\bigl|\om(x)-\om(t)\bigr|\,dx\Bigr)\,dt
\leq\frac{M}{\int_{a}^{b}\om}\int_{a}^{b}\bigl(\om(t)-\om(t-c)\bigr)\,dt
\-\\\-
=\frac{M}{\int_{a}^{b}\om}\int_{b-c}^{b}\om(t)\,dt
-\frac{M}{\int_{a}^{b}\om}\int_{a-c}^{a}\om(t)\,dt
\leq Mc\frac{\om(b)}{\int_{a}^{b}\om},
}
which tends to 0 as $b-a\ras\infty$.
\endproof
\section{S-noncom}{Noncommuting \GLSs}

If \GLSs\ $\T_{1},\ld,\T_{k}$ do not commute,
the situation becomes much more complicated.
Recall that we introduced the notions of ergodicty and joint ergodicity
of sequences of transformations
(Definitions~\rfrn{D-erg} and \rfrn{D-jerg} above)
with respect to an arbitrary fixed F{\o}lner sequence in $\Z$.
However, for commuting \GLSs\
the property of being ergodic or jointly ergodic 
has turned out to be F{\o}lner sequence independent
(see Remarks~\rfrn{R-indep} and \rfrn{R-jindep}).
An example in \brfr{BB2} shows that this is no longer the case
if $\T_{i}$ do not commute,
even in the conventional case $\T_{i}(n)=T_{i}^{n}$.
It follows that one cannot expect to have a criterion of joint ergodicity
in terms of ergodicity of a certain collection of sequences of transformations,
unless the ergodicity of these sequences is itself F{\o}lner sequence dependent.

One has nevertheless the following generalization of Theorem~2.1 from \brfr{BB2}:
\theorem{P-noncom}{}{}
Let $G_{1},\ld,G_{k}$ be several commutative groups 
of measure preserving transformations of $X$,
and for each $i=1,\ld,k$ let $\T_{i}$ be a \GLS\ in $G_{i}$.
Then $\T_{1},\ld,\T_{k}$ are jointly ergodic
iff $\T_{1},\ld,\T_{k}$ are ergodic
and $\Clim_{n}\int_{X}\prod_{i=1}^{k}\T_{i}(n)f_{i}\,d\mu=\prod_{i=1}^{k}\int_{X}f_{i}\,d\mu$
for any $f_{1},\ld,f_{k}\in L^{\infty}(X)$.
\endtheorem
\proof{}
The ``only if'' direction is clear;
we will prove the ``if'' statement.
Let $f_{1},\ld,f_{k}\in L^{\infty}(X)$, with $|f_{i}|\leq 1$ for all $i$.
First, assume that for some $i$, 
$f_{i}$ is in the $\Hwm$ space corresponding to the group $G_{i}$.
We will assume that $i=1$;
then $\T_{1}$ is weakly mixing on $f_{1}$ by \rfr{P-cwm}.
Let $\eps>0$,
and let a Bohr set $H\sle Z$, transformations $S_{i}\in G_{i}$ for $i=1,\ld,k$, 
and sets $E_{h}\sle\Z$ for $h\in H$ be as in \rfr{P-perT}(iv).
Then for any $h_{1},h_{2}\in H$,
\lequ{
\Bigl|\Clim_{n}\Blan\prod_{i=1}^{k}\T_{i}(n+h_{1})f_{i},
\prod_{i=1}^{k}\T_{i}(n+h_{2})f_{i}\Bran\Bigr|
\leq\Bigl|\Clim_{n}\int_{X}\prod_{i=1}^{k}\T_{i}(n)\bigl(\T_{i}(h_{1})S_{i}f_{i}\cd
\T_{i}(h_{2})S_{i}\af_{i}\bigr)\,d\mu\Bigr|+2\eps
\-\\\-
=\Bigl|\prod_{i=1}^{k}\int_{X}\T_{i}(h_{1})S_{i}f_{i}\cd
\T_{i}(h_{2})S_{i}\af_{i}\,d\mu\Bigr|+2\eps
\leq\Bigl|\int_{X}\T_{1}(h_{1})S_{1}f_{1}\cd\T_{1}(h_{2})S_{1}\af_{1}\,d\mu\Bigr|+2\eps.
}
Since $\Dlim_{h}\lan\T_{1}(h)S_{1}f,f'\ran=0$ for any $f'\in L^{2}(X)$,
we can construct an infinite set $B\sle H$ such that
$\bigl|\int_{X}\T_{1}(h_{1})S_{1}f_{1}\cd\T_{1}(h_{2})S_{1}\af_{1}\,d\mu\bigr|<\eps$
for any distinct $h_{1},h_{2}\in B$.
By \rfr{P-vdC}, $\Climsup_{\|\cd\|,n}\prod_{i=1}^{k}\T_{i}(n)f_{i}<\sqrt{3\eps}$.
Since $\eps$ is arbitrary, $\Clim_{n}\prod_{i=1}^{k}\T_{i}(n)f_{i}=0$.

Now assume that for each $i$, $\T_{i}$ acts on $f_{i}$ in a compact way.
We then may assume that, for each $i$, $f_{i}$ is a nonconstant eigenfunction of $G_{i}$,
and so, $\T_{i}(n)f_{i}=\lam_{i}(n)f_{i}$, $n\in\Z$,
for some \GLS\ $\lam_{i}$ in $\{z\in\C:|z|=1\}$.
In this case, $\prod_{i=1}^{k}\T_{i}(n)f_{i}=\lam(n)\prod_{i=1}^{k}f_{i}$,
where $\lam(n)=\prod_{i=1}^{k}\lam_{i}(n)$.
Since $\Clim_{n}\int_{X}\prod_{i=1}^{k}\T_{i}(n)f_{i}\,d\mu=0$,
we have $\Clim_{n}\lam(n)=0$,
and so, $\Clim_{n}\prod_{i=1}^{k}\T_{i}(n)f_{i}=0$.%
\endproof
\bibliography{}
\biblleft=13mm
\bibart Berend/B
a:B. Berend
t:Joint ergodicity and mixing
j:J. d'Analyse Math.
n:45
y:1985
p:255-284
*
\bibart BB1/BBe1
a:D. Berend and V. Bergelson
t:Jointly ergodic measure preserving transformations
j:Israel J. Math.
n:49
y:1984
p:\no{4}, 307-314
*
\bibart BB2/BBe2
a:D. Berend and V. Bergelson
t:Characterization of joint ergodicity for non-commuting transformations
j:Israel J. Math.
n:56
y:1986
p:\no{1}, 123-128
*
\bibart B-pet/Be
a:V. Bergelson
t:Weakly mixing PET
j:Ergodic Theory and Dynamical Systems
n:7
y:1987
p:\no{3}, 337-349
*
\bibook B-Gor/BeG
a:V. Bergelson and A. Gorodnik
t: Weakly mixing group actions: a brief survey and an example
i:Modern Dynamical Systems and Applications, 3-25,
Cambridge Univ. Press, New York, 2004
*
\bibart BHa/BeH 
a:V. Bergelson and I.J. H{\aa}land Knutson
t:Weak mixing implies weak mixing of higher orders along tempered functions
j:Ergodic Theory and Dynamical Systems
n:29
y:2009
p:\no{5}, 1375-1416
*
\bibart psz/BeL1
a:V. Bergelson and A. Leibman
t:Polynomial extensions of van der Waerden's and Szemeredi's theorems
j:Journal of AMS 
n:9 
y:1996
p:725-753 
*
\bibart sko/BeL2
a:V. Bergelson and A. Leibman
t:Distribution of values of bounded generalized polynomials
j:Acta Mathematica
n:198
y:2007
p:155-230
*
\bibart dcc/BeLM 
a:V. Bergelson, A. Leibman, and C.G. Moreira
t:From discrete- to continuous-time ergodic theorems
j:Ergodic Theory and Dynamical Systems
n:32
y:2012
p:\no{2}, 383-426
*
\bibook BM-PolySz/BeMc
a:V. Bergelson and R. McCutcheon
t:Uniformity in the polynomial Szemerédi theorem
i:Ergodic theory of $\Z^{d}$ actions, 273-296,
London Math. Soc. Lecture Note Ser., 228, 
Cambridge Univ. Press, Cambridge, 1996
*
\bibart B-Ros/BeR
a:V. Bergelson and J. Rosenblatt
t:Mixing actions of groups
j:Illinois J. Math
n:32
y:1988
p:\no{1}, 65-80
*
\bibart Fra/F 
a:N. Frantzikinakis
t:Multiple recurrence and convergence for Hardy sequences of polynomial growth
j:J. d'Analyse Math.
n:112
y:2010
p:79-135
*
\bibart FHK/FHoK
a:N. Frantzikinakis, B. Host, and B. Kra
t:Multiple recurrence and convergence for sequences related to prime numbers
j:J. Reine Angew. Math.
n:611
y:2007
p:131-144
*
\bibart FraKra/FK
a:N. Frantzikinakis and B. Kra
t:Ergodic averages for independent polynomials and applications
j:J. Lond. Math. Soc.
n:74
y:2006
p:131-142
*
\bibart F-Sz/Fu
a:H. Furstenberg
t:Ergodic behavior of diagonal measures and a theorem of Szemer\'edi on arithmetic progressions
j:J. d'Analyse Math.
n:31
y:1977
p:204-256
*
\bibart GT/GT
a:B. Green and T. Tao
t:Linear equations in primes
j:Ann. of Math. (2)
n:171
y:2010
p:\no{3}, 1753-1850
*
\bibart GTZ/GTZ
a:B. Green, T. Tao, and T. Ziegler
t:An inverse theorem for the Gowers $U^{s+1}$-norm
j:Ann. of Math.
n:176
y:2012
p:\no(2), 1231-1372
*
\bibart vN-K/vNK
a:J. von Neumann and B.O. Koopman
t:Dynamical systems of continuous spectra
j:Proc. Nat. Acad. Sci. 
n:18 
y:1932
p:255-263
*
\bibartx Sun/S
a:W. Sun
t:Multiple recurrence and convergence for certain averages along shifted primes
x:Available at arXiv:1303.3902
*
\endbibliography
\finish
\bye

%% file: xpehfnt.tex
\def\fnt#1{\csname !f-#1\endcsname}
\def\deffnt#1{\expandafter\font\csname !f-#1\endcsname}
\def\scal#1{ scaled \magstep#1}
\deffnt{rm}=cmr10
\deffnt{rm.1}=cmr10 \scal1
\deffnt{rm.2}=cmr10 \scal2
\deffnt{rm9}=cmr9
\deffnt{rm8}=cmr8
\deffnt{rm7}=cmr7
\deffnt{rm5}=cmr5
\deffnt{bf}=cmbx10
\deffnt{bf.1}=cmbx10 \scal1
\deffnt{bf.2}=cmbx10 \scal2
\deffnt{bf9}=cmbx9
\deffnt{bf7}=cmbx7
\deffnt{bf5}=cmbx5
\deffnt{it}=cmti10
\deffnt{it.1}=cmti10 \scal1
\deffnt{it9}=cmti9
\deffnt{it7}=cmti7
\deffnt{mi}=cmmi10
\deffnt{mi.1}=cmmi10 \scal1
\deffnt{mi9}=cmmi9
\deffnt{mi7}=cmmi7
\deffnt{mi5}=cmmi5
\deffnt{sy}=cmsy10
\deffnt{sy9}=cmsy9
\deffnt{ex}=cmex10
\deffnt{ex9}=cmex9
\deffnt{sl}=cmsl10
\deffnt{sl8}=cmsl8
\deffnt{sf}=cmss10
\deffnt{sf.2}=cmss10 \scal2
\deffnt{sc}=cmcsc10
\deffnt{sc9}=cmcsc9
\deffnt{sc.1}=cmcsc10 \scal1
\deffnt{mb}=msbm10
\deffnt{mb7}=msbm7
\deffnt{mb5}=msbm5
\newfam\mbfam
\textfont\mbfam=\fnt{mb}
\scriptfont\mbfam=\fnt{mb7}
\scriptscriptfont\mbfam=\fnt{mb5}
\def\mb{\fam\mbfam\fnt{mb}}
\newfam\lgfam
\textfont\lgfam=\fnt{mi.1}
\scriptfont\lgfam=\fnt{mi9}
\scriptscriptfont\lgfam=\fnt{mi7}

\newfam\smfam
\textfont\smfam=\fnt{mi9}
\scriptfont\smfam=\fnt{mi7}
\scriptscriptfont\smfam=\fnt{mi5}
\deffnt{gt}=eufm10
\deffnt{gt5}=eufm5

%% file: xpehcmd.tex
\def\cmd#1{\csname #1\endcsname}
\def\defcmd#1{\expandafter\def\csname #1\endcsname}
\def\edefcmd#1{\expandafter\edef\csname #1\endcsname}
\def\ifundefined#1{\expandafter\ifx\csname #1\endcsname\relax}
\edef\warning#1{{\newlinechar=`|\message{|*********** #1}}}
\def\npar{\relax\par\noindent}
\newbox\xbox
\def\ifhempty#1#2#3{\setbox\xbox=\hbox{#1}\ifdim\wd\xbox=0pt\relax
                                                           #2\else #3\fi}
\def\lph#1{\hbox to 0pt{\hss #1}}
\def\rph#1{\hbox to 0pt{#1\hss}}
\def\ulph#1{\vbox to 0pt{\vss\hbox to 0pt{\hss #1}}}
\def\urph#1{\vbox to 0pt{\vss\hbox to 0pt{#1\hss}}}
\def\dlph#1{\vbox to 0pt{\hbox to 0pt{\hss #1}\vss}}
\def\drph#1{\vbox to 0pt{\hbox to 0pt{#1\hss}\vss}}
\def\nasp{\spacefactor=1000}
\def\small{\def\rm{\fnt{rm9}\fam\smfam}\def\it{\fnt{it9}}\def\bf{\fnt{bf9}}%
\textfont1=\fnt{mi9}\textfont2=\fnt{sy9}\textfont3=\fnt{ex9}%
\baselineskip=.8\baselineskip\parindent=.8\parindent%
\abovedisplayskip=.6\abovedisplayskip\belowdisplayskip=.6\belowdisplayskip%
\stateskip=.6\stateskip
\rm}
\def\nohyphen{\hyphenpenalty=5000 \exhyphenpenalty=5000
 \tolerance=5000 \pretolerance=5000}
\def\center{\leftskip=0pt plus 1fill\relax \rightskip=\leftskip\relax
 \parfillskip=0pt\relax \parindent=0pt\relax \nohyphen}
\def\\{\ifhmode\hfil\break\else\ifvmode\vskip\baselineskip\fi\fi}
\def\-{\hfill}
\def\vsp#1{\vadjust{\kern #1}}

\def\vbreak#1#2{\par\ifdim\lastskip<#2\removelastskip\penalty#1\vskip#2\fi}
\def\vbreakn#1#2{\vbreak{#1}{#2}\npar}

\newcount\multilabels \newcount\undeflabels
\newread\filelabelsold \newwrite\filelabels
\edef\filelabelsname{\jobname.lbl}
\def\getlabels{\openin\filelabelsold=\filelabelsname\relax
 \ifeof\filelabelsold 
  \warning{The file of labels \filelabelsname\ does not exist}
  \else\closein\filelabelsold
   {\catcode`/=0 \globaldefs=1 \input \filelabelsname\relax}\fi}
\def\openlabels{\immediate\openout\filelabels=\filelabelsname\relax}
\def\closelabels{\closeout\filelabels}

\def\label#1#2#3{
\ifx\printlabel Y\lnote{\fnt{rm7}#1}\fi%
\immediate\write\filelabels{%
/defcmd{!tl-#1}{#2}/defcmd{!nl-#1}{#3}}%
\ifundefined{!dl-#1}\defcmd{!dl-#1}{!earlier-defined!}%
\else\warning{Label #1 is multiply defined.}\advance\multilabels by 1\fi%
\edefcmd{!tl-#1}{#2}\edefcmd{!nl-#1}{#3}}
\def\xlabel#1#2#3#4{
\ifx\printlabel Y\ldnote{\fnt{rm7}#1}\fi%
\immediate\write\filelabels{%
/defcmd{!tl-#1}{#2}/defcmd{!nl-#1}{#3}/defcmd{!xl-#1}{#4}}%
\ifundefined{!dl-#1}\defcmd{!dl-#1}{!earlier-defined!}%
\else \warning{LABEL #1 is multiply defined.} \advance\multilabels by 1\fi%
\edefcmd{!tl-#1}{#2}\edefcmd{!nl-#1}{#3}\edefcmd{!xl-#1}{#4}}
\def\undelabel#1{\warning{Label #1 has not been defined.}%
\global\advance\undeflabels by 1{?#1?}}
\def\rfrn#1{\ifundefined{!nl-#1}\undelabel{#1}\else\cmd{!nl-#1}\fi}

\def\rfrx#1{\ifundefined{!xl-#1}\undelabel{#1}%
 \else\cmd{!xl-#1}\fi}
\def\rfr#1{\ifundefined{!tl-#1}\undelabel{#1}%
 \else\cmd{!tl-#1}~\cmd{!nl-#1}\fi}

\def\title#1{\vbreak{-200}{1cm}{\center\fnt{bf.1}#1\par}}
\def\author#1{\vbreak{100}{5mm}{\center\fnt{rm}#1\par}}
\def\titledate{\vbreak{100}{5mm}\centerline{\fnt{it}\date}}

\newskip\absskip
\absskip=20mm
\long\def\abstract#1{\vbreak{50}{5mm}\centerline{\fnt{bf}Abstract}
 \kern3mm{\leftskip=\absskip\rightskip=\leftskip
 \small #1\vbreak{-9000}{4mm}}}

\def\support#1{\footnote{}{{\small#1\hfil}}}

\ifx\headdate Y{\headline={{\fnt{rm7}\date}\hfil}}\fi
\def\date{\ifcase\month%
\or January\or February\or March\or April\or May\or June\or July%
\or August\or September\or October\or November\or December\fi%
~\number\day,~\number\year}
\def\address#1#2#3{\vbreakn{-9000}{1mm}#1, #2\ifhempty{#3}{}{\\{\it e-mail:} #3}\par}
\def\myaddress{\rm%
\setbox\xbox=\hbox{Department of mathematics}%
\vtop{\hsize=\wd\xbox\parindent=0pt
Department of Mathematics\\
The Ohio State University\\
Columbus, OH 43210, USA\\
\hbox to \wd\xbox{{\it e-mail}: leibman@math.ohio-state.edu\hss}}}


\newcount\fnotenum
\def\fnote#1{\global\advance\fnotenum by 1
 {\parindent=10mm \parfillskip=0pt plus 1fill\relax
  \footnote{$^{(\the\fnotenum)}$}{{\small #1}}}}

\def\lnote#1{\ifvmode\ulph{#1 }\else\vadjust{\ulph{#1 }}\fi}
\def\ldnote#1{\ifvmode\ulph{#1 }\else\vadjust{\dlph{#1 }}\fi}
\newcount\sectionnum \sectionnum=-1 
\newcount\subsectionnum \subsectionnum=0 
\newcount\subsubsectionnum \subsubsectionnum=0
\newskip\sectionskip \sectionskip=6mm
\newskip\aftersectionskip \aftersectionskip=3mm
\newcount\sectionpenalty \sectionpenalty=-9000
\def\xsection#1{\global\advance\sectionnum by 1%
\global\subsectionnum=0\global\subsubsectionnum=0%
\ifx\nonewsecequation Y\else\global\equationnum=0\fi%
\ifx\nonewsecstate Y\else\global\statenum=0\fi%
\label{#1}{Section}{\the\sectionnum}}
\def\section#1#2{%
\vbreak{\sectionpenalty}{\sectionskip}\xsection{#1}%
\vbox{\ifx\nosectioncenter Y\npar\else\center\fi\fnt{bf}
\ifx\nosectionnum Y\else\rfrn{#1}.{\nasp} \fi#2\par}%
\vbreak{10000}{\aftersectionskip}}
\newskip\subsectionskip \subsectionskip=3mm
\newcount\subsectionpenalty \subsectionpenalty=-7000
\def\thesubsectionnum{\ifx\nosection Y\the\subsectionnum
\else\the\sectionnum.\the\subsectionnum\fi}
\def\xsubsection#1{\global\advance\subsectionnum by 1%
\global\subsubsectionnum=0%
\label{#1}{subsection}{\thesubsectionnum}}
\def\subsection#1{%
\vbreak{\subsectionpenalty}{\subsectionskip}
\ifx\nosubsectionnum Y\par\else
\noindent\xsubsection{#1}{\fnt{bf}\rfrn{#1}.{\nasp}}\fi}
\def\thesubsubsectionnum{\thesubsectionnum.\the\subsubsectionnum}
\def\xsubsubsection#1{\global\advance\subsubsectionnum by 1%
\label{#1}{subsection}{\thesubsubsectionnum}}
\def\subsubsection#1{%
\vbreak{-5000}{1.5mm}\noindent\xsubsubsection{#1}{\fnt{bf}\rfrn{#1}.{\nasp}}}
\newcount\equationnum   
\def\theequation{\ifx\nosection Y\the\equationnum
\else\the\sectionnum.\the\equationnum\fi}
\def\equn#1#2{\global\advance\equationnum by 1%
\label{#1}{Equation}{\theequation}%
$$\vcenter{\equalign{#2}}\eqno(\theequation)$$}
\def\equ#1{$$\vcenter{\equalign{#1}}$$}
\def\equalign#1{\let\\=\cr\let\-=\hfill
\ialign{&\hfil$\dsp ##$\hfil\cr#1\crcr}}
\def\lequ#1{$$\vcenter{\lequalign{#1}}$$}
\def\lequalign#1{\let\\=\cr\let\-=\hfill
\ialign{&\hbox to \hsize{$\dsp ##$\hfil}\cr#1\crcr}}
\def\matalign#1#2{{\let\\=\cr\let\-=\hfill%
\ialign{\hfil$##$\hfil&&#1\hfil$##$\hfil\cr#2\crcr}}}

\def\frfr#1{\hbox{(\rfrn{#1})}}
\newskip\stateskip\stateskip=2mm
\newcount\statenum\statenum=0
\def\thestate{\ifx\afterstate Y\the\statenum\else\thesubsectionnum\fi}
\def\state#1#2#3#4{\global\advance\statenum by 1\xlabel{#2}{#1}{\thestate}{#3}%
{\bf #1\ifhempty{#3}{\ifx\afterstate Y\ \thestate\fi}{ \box\xbox}.}{\nasp} #4%
\begingroup}
\def\nstate#1#2{{\bf #1\ifhempty{#2}{}{ \box\xbox}.}{\nasp}\begingroup}
\def\endstate{\par\endgroup\vbreak{-7000}{\stateskip}}
\def\ltheorem#1#2#3{\state{Theorem}{#1}{#2}{#3}\it}
\def\theorem#1#2#3{\vbreakn{-7000}{\stateskip}\ltheorem{#1}{#2}{#3}}
\def\endtheorem{\endstate}
\def\llemma#1#2#3{\state{Lemma}{#1}{#2}{#3}\it}
\def\lemma#1#2#3{\vbreakn{-7000}{\stateskip}\llemma{#1}{#2}{#3}}
\def\endlemma{\endstate}
\def\lproposition#1#2#3{\state{Proposition}{#1}{#2}{#3}\it}
\def\proposition#1#2#3{\vbreakn{-7000}{\stateskip}\lproposition{#1}{#2}{#3}}
\def\endproposition{\endstate}
\def\lcorollary#1#2#3{\state{Corollary}{#1}{#2}{#3}\it}
\def\corollary#1#2#3{\vbreakn{-7000}{\stateskip}\lcorollary{#1}{#2}{#3}}
\def\endcorollary{\endstate}


\def\lremark#1{\nstate{Remark}{#1}}
\def\remark#1{\vbreakn{-7000}{\stateskip}\lremark{#1}}
\def\endremark{\endstate}
\def\lnremark#1#2{\state{Remark}{#1}{#2}{}}
\def\nremark#1#2{\vbreakn{-7000}{\stateskip}\lnremark{#1}{#2}}
\def\endnremark{\endstate}

\def\ldefinition#1{\nstate{Definition}{#1}}
\def\definition#1{\vbreakn{-7000}{\stateskip}\ldefinition{#1}}
\def\enddefinition{\endstate}
\def\lndefinition#1#2{\state{Definition}{#1}{#2}{}}
\def\ndefinition#1#2{\vbreakn{-7000}{\stateskip}\lndefinition{#1}{#2}}
\def\endndefinition{\endstate}
\def\lexample#1{\nstate{Example}{#1}}
\def\example#1{\vbreakn{-7000}{\stateskip}\lexample{#1}}
\def\endexample{\endstate}

\def\lproof#1{\nstate{Proof}{#1}}
\def\proof#1{\vbreakn{-7000}{\stateskip}\lproof{#1}}

\def\endproof{\endprr\endstate}
\def\enprule{\vrule height1mm depth1mm width2mm}
\def\endprr{\discretionary{}{\kern\hsize}{\kern 3ex}\llap{\enprule}}
\def\frgdsp{\par\penalty10000
 \vskip-\belowdisplayskip\kern-2mm\noindent\hbox to \hsize{\hfil}}
\newcount\biblnum \newcount\biblpenalty
\def\bibliography#1{\vbreak{-9000}{6mm}
\biblpenalty=10000\biblnum=0%
\line{\fnt{bf}\ifhempty{#1}{Bibliography}{\box\xbox}\hfil}\par\kern2mm
\begingroup\parskip=0pt\parindent=0pt\frenchspacing}
\def\endbibliography{\par\endgroup}
\newskip\biblleft \biblleft=10mm
\def\bibxitem#1/#2(#3) #4{\advance\biblnum by 1%
\vbreak{\biblpenalty}{.8mm}\biblpenalty=-9000%
\hangindent=\biblleft \hangafter=1
\noindent{\small\rlap{\brfr{#1}}\hskip\hangindent
\xlabel{#1}{#3}{\the\biblnum}{#2}#4}}
\def\bibitem#1/#2(#3) #4#5{\bibxitem#1/#2(#3) {{\frenchspacing #4}, #5}}
\def\bibx#1/#2 {\advance\biblnum by 1{#1/#2}
\xlabel{#1}{text}{\the\biblnum}{#2}}
\def\bibart#1/#2 a:#3 t:#4 j:#5 n:#6 y:#7 p:#8 *{%
\bibitem{#1}/{#2}(paper) {#3}{{#4}, {\it #5\/}
   {\bf #6} (#7), #8.}}
\def\bibartp#1/#2 a:#3 t:#4 *{%
\bibitem{#1}/{#2}(paper) {#3}{{#4}, {in preparation}.}}
\def\bibarts#1/#2 a:#3 t:#4 *{%
\bibitem{#1}/{#2}(paper) {#3}{{#4}, {submitted}.}}
\def\bibarta#1/#2 a:#3 t:#4 j:#5 *{%
\bibitem{#1}/{#2}(paper) {#3}{{#4}, {\rm to appear in \it #5}.}}
\def\bibartn#1/#2 a:#3 t:#4 j:#5 *{%
\bibitem{#1}/{#2}(paper) {#3}{{#4}, #5.}}
\def\bibartpr#1/#2 a:#3 t:#4 *{%
\bibitem{#1}/{#2}(paper) {#3}{{#4}, {preprint}.}}
\def\bibartx#1/#2 a:#3 t:#4 x:#5 *{%
\bibitem{#1}/{#2}(paper) {#3}{{#4}, #5.}}
\def\bibook#1/#2 a:#3 t:#4 i:#5 *{%
\bibitem{#1}/{#2}(book) {#3}{{\it #4\/}, #5.}}
\def\no#1{no.{\nasp}~{#1}}

\def\brfr#1{\hbox{\rm[\ifx\biblNUM Y\rfrn{#1}\else\rfrx{#1}\fi]}}
\def\start{
 \warning{|****** Start ******|}
 \getlabels\openlabels\multilabels=0\undeflabels=0
 \null}
\def\finish{\par\closelabels
 \warning{|****** Finish ******}
 \ifnum\multilabels>0\warning{\the\multilabels_multidefined labels!}\fi
 \ifnum\undeflabels>0\warning{\the\undeflabels_undefined labels!}
                     \warning{*************************}
                     \warning{Try to run TeX once again!}
                     \warning{*************************}\fi}
\output={\shipout\vbox{\makeheadline\pagebody\makefootline}
\ifx\doublepage Y\shipout\vbox{}\fi
 \advancepageno
 \ifnum\outputpenalty>-20000\else\dosupereject\fi}
\def\sdup#1#2{{\scr #1 \atop \scr #2}}

\def\Gal{\lnote{{\bf V\kern 1cm\relax}}}
\def\comment#1{\ifhmode\\\fi
\hbox to \hsize{\hss\vbox{\advance\hsize by 10mm
\baselineskip=.7\baselineskip\npar\ulph{\fnt{bf.2}V}\fnt{bf7}#1}}}

\def\rest#1{\raise-2pt\hbox{$|_{#1}$}}

\def\notdvd{\mathrel{\setbox\xbox=\hbox{$\big|$}%
\hbox to 0pt{\hbox to\wd\xbox{\hss/\hss}\hss}\big|}}

\def\frac#1#2{{#1\over #2}}
\def\matr#1{{\let\\=\cr\left(\matrix{#1\crcr}\right)}}
\def\vect#1{{\let\\=\cr\left(\matrix{#1\crcr}\right)}}
\def\smatr#1{{\baselineskip=2pt\lineskip=2pt
\left(\vcenter{\let\\=\cr\let\-=\hfill
\ialign{\hfil$\scr##$\hfil&&\hfil\kern2pt$\scr##$\hfil\cr#1\crcr}}\right)}}

\def\svd{\vbox to 2.4mm{\nolineskip
\kern.3pt\hbox{.}\vfil\hbox{.}\vfil\hbox{.}\kern.3pt}}
\def\vd{\vbox to 3.2mm{\nolineskip
\kern1pt\hbox{.}\vfil\hbox{.}\vfil\hbox{.}\kern1pt}}
\def\lvd{\vbox to 7mm{\nolineskip
\kern1pt\hbox{.}\vfil\hbox{.}\vfil\hbox{.}\vfil\hbox{.}
\vfil\hbox{.}\vfil\hbox{.}\kern1pt}}
\def\rvd{\vbox to 2.4mm{\nolineskip
\kern.3pt\hbox{.}\vfil\hbox{\kern3.5pt.}\vfil\hbox{\kern7pt.}\kern.3pt}}
\def\lld{\hbox to 5mm{\kern2pt.\hfil.\hfil.\kern2pt}}

\def\comp{\mathord{\hbox{$\scr\circ$}}}
\def\semprod{\mathrel{\times\kern-2pt
\vbox{\hrule width.5pt height4pt depth0pt\kern.5pt}\kern1pt}}

\def\vkeep#1#2{\vrule width 0pt height #1 depth #2}
\def\vs#1{\vkeep{0pt}{#1}}
\def\dsc{\discretionary{}{}{}}
\long\def\omit#1\endomit{\par\vbox{%
\hrule\vskip2mm\hfil\vdots\hfil\vskip2mm\hrule}}
\long\def\ignore#1\endignore{}
\def\R{{\mb R}}
\def\C{{\mb C}}
\def\Z{{\mb Z}}
\def\N{{\mb N}}
\def\Q{{\mb Q}}
\def\PP{\raise2.2pt\hbox{\fnt{mi.1}\char"7D}}

\def\mod{\mathop{\rm mod}}

\def\const{\mathop{\hbox{\rm const}}}


\def\nolineskip{\baselineskip=0pt\lineskip=0pt}

\let\dsp=\displaystyle
\let\txt=\textstyle
\let\ovr=\overline

\let\scr=\scriptstyle

\let\sln=\subset
\let\sle=\subseteq

\let\sm=\setminus

\let\col=\colon
\let\lan=\langle
\let\ran=\rangle

\let\ld=\ldots
\let\vd=\vdots
\let\cd=\cdot

\let\ra=\longrightarrow
\let\ras=\rightarrow

\let\tens=\otimes
\let\alf=\alpha
\let\bet=\beta
\let\Gam=\Gamma

\let\del=\delta

\let\phi=\varphi
\let\sig=\sigma
\let\eps=\varepsilon
\let\lam=\lambda

\let\om=\omega
